\documentclass{article}
   \usepackage{amssymb}  
 \usepackage[totalwidth=410pt, totalheight=620pt]{geometry}
\newcommand{\eps}{\epsilon}

\newcommand{\La}{\Lambda}

\newcommand{\Th}{\Theta}

\newcommand{\mb}{\mbox}

\newcommand{\Mklz}[2]{\left\{\left.\;#1\;\right|\; #2\;\right\}}
\newcommand{\F}{\mathbb{F}}

\newcommand{\N}{\mathbb{N}}
\newcommand{\Nn}{\mathbb{N}_0}

\newcommand{\Z}{\mathbb{Z}}
\newcommand{\We}{\mbox{$\cal W$}}

\newcommand{\GD}{\mbox{$\widehat{G}$}}

\newcommand{\Proof}{\mbox{\bf Proof: }}
\newcommand{\End}{\mb{}\hfill\mb{$\square$}\\}

\newcommand{\ti}{\tilde}

\newcommand{\idem}{\,e\,}
\begin{document}
\newtheorem{Theorem}{Theorem}[section]
\newtheorem{Def}[Theorem]{Definition}
\newtheorem{Prop}[Theorem]{Proposition}
\newtheorem{Prop+Def}[Theorem]{Proposition+Definition}
\newtheorem{Cor}[Theorem]{Corollary}
\newtheorem{Rem}[Theorem]{Remark}
\newtheorem{Rems}[Theorem]{Remarks}
\newtheorem{Lemma}[Theorem]{Lemma}
%
%
\title{The maximal chains of the extended Bruhat orders on the ${\cal W}\times {\cal W}$-orbits of an 
infinite Renner monoid}
\author{Claus Mokler\thanks{Supported by the Deutsche Forschungsgemeinschaft}\\\\ Universit\"at Wuppertal, Fachbereich Mathematik\\  Gau\ss stra\ss e 20\\ D-42097 Wuppertal, Germany\vspace*{1ex}\\ 
          mokler@math.uni-wuppertal.de}
\date{}
\maketitle
\begin{abstract}\noindent
Let $(\We,S)$ be a Coxeter system. For $\eps,\delta\in\{+,-\}$ we introduce and investigate combinatorially certain 
partial orders $\leq_{\eps\delta}$, called extended Bruhat orders, on a $\We\times\We$-set $\We(N,C)$, which depends on 
$\We$, a subset $N\subseteq S$, and a component $C\subseteq N$. We determine the length of the maximal chains 
between two elements $x,y\in\We(N,C)$, $x\leq_{\eps\delta}y$.\\
These posets generalize $\We$ equipped with its Bruhat order. They include the $\We\times\We$-orbits of the 
Renner monoids of reductive algebraic monoids and of some infinite dimensional generalizations which are equipped with the 
partial orders obtained by the closure relations of the Bruhat and Birkhoff cells. 
They also include the $\We\times\We$-orbits of certain posets obtained by generalizing the closure relation of the Bruhat 
cells of the wonderful compactification. 
\end{abstract}
\vspace*{1ex}
\section*{Introduction}
Of particular importance in the theory of reductive algebraic monoids, which has been developed mainly by 
M. S. Putcha and L. E. Renner, is the Bruhat decomposition and its associated structures, \cite{Pu1}, \cite{Pu2}, 
\cite{Re1}, \cite{Re2}:
Let $M$ be a reductive algebraic monoid. Let $G$ be its reductive unit group. Let $T$ be a maximal torus of $G$ and 
$B=B^+$, $B^-$ be opposite Borel subgroups containing $T$. Equip $M$ with its natural $G\times G$-action. The Bruhat 
decomposition is the $B\times B$-orbit decomposition
\begin{eqnarray}\label{M=BB}
   M=\dot{\bigcup_{x\in {\cal R}}} BxB 
\end{eqnarray}
where ${\cal R}=\overline{N_G(T)}/T$ is the Renner monoid. The Renner monoid contains the Weyl group $\We=N_G(T)/T$ 
as unit group. Equipped with its natural $\We\times\We$-action, the $\We\times\We$-orbits can be parametrized by a 
certain lattice of idempotents $\Lambda\subseteq {\cal R}$, called a cross section lattice, i.e.,
\begin{eqnarray}\label{R=WW}
 {\cal R}\;=\;\dot{\bigcup_{e\in\Lambda}}\We e\We\;.
\end{eqnarray} 
To this decomposition correspond the decompositions
\begin{eqnarray}\label{M=GG}
      M= \dot{\bigcup_{e\in\Lambda}} G e G \quad\mb{ and }\quad GeG=\dot{\bigcup_{x\in {\cal W} e {\cal W}}} BxB  \;.
\end{eqnarray}
A classification of the possible cross section lattices has been achieved only in special cases, \cite{PuRe}.
On the other hand a $\We\times\We$-orbit $\We e\We$, $e\in\La$, is easy to describe. It is determined by $\We$, and by 
the normalizer $N(e) :=\Mklz{w\in\We}{we=ew}$ and the centralizer $C(e) := \Mklz{w\in\We}{we = e = ew  }$ of $e$. Both 
are standard parabolic subgroups of $\We$, i.e., $N(e) :=\We_N$ and $C(e) := \We_C$, and $C$ is a component of $N$.\\
The closure relation of the $G\times G$-orbits in the decomposition (\ref{M=GG}) of $M$, transfered to $\La$, is 
given by the partial order of the cross section lattice $\La$.
The closure relation of the Bruhat cells of the Bruhat decomposition (\ref{M=BB}) of $M$, transfered to the 
Renner monoid ${\cal R}$, is called Bruhat-Chevalley order. It has been investigated in a series of papers. 
L. E. Renner showed in \cite{Re1} by an algebraic geometric proof that all maximal chains between two elements 
$x,y\in{\cal R}$, $x\leq y$, have the same length. He introduced and investigated in \cite{Re2} a natural, algebraic 
geometrically defined length function on ${\cal R}$. 
E. A. Pennel, M. S. Putcha, and L. E. Renner obtained in \cite{PePuRe} an algebraic description of the 
Bruhat-Chevalley order and the length function. M. S. Putcha investigated in \cite{Pu3} the lexicographic shellability 
and the M\"obius function of the Bruhat-Chevalley order restricted to the $\We\times\We$-orbits $\We e\We$. In 
particular he showed that in the case of $C(e)= 1 $ the restricted Bruhat-Chevalley order is CL-shellable 
and Eulerian.\vspace*{1ex}\\ 
Let $G$ be a semisimple algebraic group of adjoint type. L. E. Renner gave in \cite{Re3} a monoid approach to the 
wonderful compactification $Wcp$ of 
\begin{eqnarray*}
       G\times G \;/ \;\mb{Diagonal of } G\times G\;.
\end{eqnarray*}
This makes the results obtained for reductive algebraic monoids available for this compactification. In particular 
there are decompositions
\begin{eqnarray*}
      Wcp = \dot{\bigcup_{x\in {\cal R}\setminus\{0\}}} BxB \quad,\quad  
      Wcp= \dot{\bigcup_{e\in\Lambda\setminus \{0\}}} G e G \quad\mb{ and }\quad 
                                          GeG=\dot{\bigcup_{x\in {\cal W} e {\cal W}}} BxB  \;,
\end{eqnarray*}
where ${\cal R}$ is a certain Renner monoid with cross section $\La$. In detail $\La\setminus\{0\}$ is isomorphic to 
the lattice of subsets of the set of simple reflections $S$ of the Weyl group $\We$. If $e(I)\in\La$ is the idempotent 
corresponding to $I\subseteq S$, then $C(e(I))=1$ and $N(e(I))=\We_I$.\vspace*{1ex}\\
At the same time T. A. Springer investigated in \cite{Sp} the intersection cohomology of the $B\times B$-orbit closures of the wonderful compactification of above. For this he determined the $B\times B$-orbits and their closure relation, 
which he called Bruhat order, explicitely by a different approach. (The poset $V$ of \cite{Sp} identifies with 
${\cal R}\setminus \{0\}$ equipped with the Bruhat-Chevalley order by mapping $[I,a,b]\in V$ to 
$b e_I a^{-1}\in {\cal R}\setminus \{0\}$. To see this use the results of the following Section \ref{mult}.) He also 
introduced a compatible length function.
He showed that most of the structures obtained in his investigation of the intersection cohomology can be 
extended combinatorially to arbitrary Coxeter groups. It remained open if a certain map $\Delta$, which generalizes a 
map related to the Verdier duality, is involutive. This property is equivalent to the existence of certain 
analogues of Kazhdan-Lusztig polynomials.\vspace*{1ex}\\ 
This question has been solved by Y. Chen and M. J. Dyer in \cite{ChDy}. In the series of papers \cite{Dy2}, 
\cite{Dy2}, and \cite{Dy3} M. J. Dyer introduced and investigated combinatorially certain generalizations of the 
Bruhat order and length function on Coxeter groups, called twisted Bruhat orders and twisted length functions.
The main aid of Y. Chen and M. J. Dyer in \cite{ChDy} is a $\We\times\We$-equivariant order embedding of the set 
${\cal R }\setminus\{ 0 \}$ together with its Bruhat order (defined for a Coxeter group) into a non-canonically 
associated Coxeter group equipped with a certain twisted Bruhat order, preserving the corresponding 
length functions up to a additive constant, and preserving the corresponding analogues of Kazhdan-Lustzig 
R-polynomials.
Y. Chen and M. J. Dyer also used this isomorphism to transfer properties of the twisted Bruhat-Chevalley order to 
the Bruhat order of ${\cal R}\setminus \{0\}$ (defined for a Coxeter group). In particular the maximal chains between 
$x,y\in {\cal R}\setminus \{0\}$, $x\leq y$, have the same length, given by the difference of the length functions of 
$x$ and $y$. They also obtained the pure EL-shellability of closed intervals of the whole ${\cal R}\setminus\{0\}$.
The proof that this map is actually an order isomorphism uses the analogues of the Kazhdan Lustzig R-polynomials for 
the twisted Bruhat-Chevalley order, Springers analogues of Kazhdan-Lustzig R-polynomials, and the properties of these 
polynomials, as well as properties of Springers function $\Delta$, requiring the whole construction of these 
things.\vspace*{1ex}\\ 
The author investigated in \cite{M1} an analogue of a reductive algebraic monoid $\widehat{G}$, whose unit group is a 
Kac-Moody group $G$. Its coordinate ring restricted to $G$ is the algebra of strongly regular functions of V. Kac and 
D. Peterson, \cite{KP}.
This monoid is a purely infinite dimensional phenomenon. In the classical case it reduces to the group $G$ itself. 
For its history please compare the introduction of \cite{M1}.\\
The monoid $\GD$ has similar structural properties as a reductive algebraic monoid. In particular there are Bruhat 
and Birkhoff decompositions. The corresponding Renner monoid $\widehat{\cal W}$, called Weyl monoid in \cite{M1}, 
is infinite in the non-classical case. It is described as follows: The cross section lattice $\La$ can be identified 
with the subsets $\Th$ of the simple reflections of $\We$, such that either $\Th$ is empty or its Coxeter diagram 
contains no component of finite type. For $e(\Th)\in\La$ corresponding to the set $\Th$ we have $C(e(\Th))=\We_\Th$ 
and $N(e(\Th))=\We_{\Th\cup\Th^\bot}$, where $\Th^\bot$ consists of the set of simple reflections which commute with 
every simple reflection of $\Th$. In difference to the cross section lattices of reductive algebraic monoids the cross 
section lattice here may contain maximal chains of different length.\\
In \cite{M3} the closure relations $\leq_{\eps\delta}$ of the Bruhat and Birkhoff cells of the decompositions 
\begin{eqnarray*}
  \GD=\bigcup_{ x\in \widehat{\cal W}} B^\eps x B^\delta \quad\mb{ where }\quad (\eps,\delta) =(+,+), (-,-),(-,+)\;,
\end{eqnarray*}
called extended Bruhat orders, have been determined. The results are similar to the case of a reductive 
algebraic monoid. The proofs are different, because most of the theorems of algebraic geometry which are used to 
investigate algebraic groups and monoids break down for these infinite dimensional varieties.
In particular the proof of L. E. Renner that all maximal chains of the Bruhat-Chevalley order between between two 
elements of the Renner monoid have the same length can not be generalized to this situation. Also there is no 
longest element of the Weyl group, which has as a consequence that the extended Bruhat order $\le_{-+}$ is 
quite different from the extended Bruhat orders $\le_{++}$, $\leq_{--}$.\\ 
The length of the maximal chains of the extended Bruhat orders between two elements of a $\We\times\We$-orbit 
$\We e(\Th)\We$ is important for determining the Krull codimension between certain Bruhat and Birkhoff cells 
contained in a $G\times G$-orbit of $\GD$. It is to expect that it will also be important for the investigation 
of a completion of the flag variety of Kashiwara \cite{Kas} resp. Pickrell \cite{Pi}, please compare the 
introduction of \cite{M2}. In this paper we determine the length of these chains combinatorially.\\
We do it in a general setting including all the $\We\times\We$-orbits of the Renner monoids equipped with their 
Bruhat(-Chevalley) order of above, obtaining a direct combinatorial proof in these cases.
In Section \ref{NP} we introduce our notation on Coxeter systems $(\We,S)$ and state some theorems on the Bruhat 
order $\leq$ of $\We$ which we use very often.
In Section \ref{mult} we investigate to which extend it is possible to multiply $a\leq b$ by $w$, $a,b,w\in\We$. The 
results of this section are used in many proofs of the following Sections \ref{BL} and \ref{MC}.
Starting with a Coxeter system $(\We,S)$, a subset $N\subseteq S$, and a component $C$ of $N$ we define in Section 
\ref{BL} for $\eps,\delta=\{+,-\}$ a relation $\leq_{\eps\delta}$ and a function $l_{\eps\delta}$ on a certain set 
$\We(N,C)$. We show that $\leq_{\eps\delta}$ is a partial order compatible with $l_{\eps\delta}$. We call these 
partial orders and functions the extended Bruhat orders and extended length functions. We investigate the extended 
Bruhat orders, in particular we give different characterizations.
In Section \ref{MC} we show that all maximal chains between two elements $x,y\in\We(N,C)$, $x\leq_{\eps\delta}y$, have 
the same length $l_{\eps\delta}(y)-l_{\eps\delta}(x)$. This also leads to the Z-Lemma for the extended Bruhat orders. 
Furthermore we obtain easy systems of relations generating the extended Bruhat orders, generalizing the system of 
generators used in the usual definition of the Bruhat order on the Weyl group. 
%
%
%
%
%
\section{Notation and Preliminaries\label{NP}}
We first introduce our notation on Coxeter systems. For the definitions compare the book \cite{Hu} 
of J. E. Humphreys:\vspace*{1ex}\\
In the whole paper $(\We,S)$ is a {\it Coxeter system} with {\it Coxeter group} $\We$ and finite set of 
{\it simple reflections} $S$. We denote by $T$ its set of {\it reflections}.\vspace*{1ex}\\
We denote by $l:\We\to \Nn$ the {\it length function} of the Coxeter system. We denote by 
$\leq$ the {\it Bruhat order} on $\We$.\vspace*{1ex}\\ 
Let $J\subseteq S$. Then $\We_J$ is the parabolic subgroup generated by $J$. $\We^J$ denotes the set of 
{\it minimal coset representatives} of $\We/\We_J$. If $w\in\We$ then $w=w^J w_J$ is the unique 
decomposition with $w^J\in\We^J$, $w_J\in\We_J$.
Similarly $\mb{}^J\We$ denotes the set of {\it minimal coset representatives} of $\We_J\backslash \We$. If 
$w\in\We$ then $w=w_J \mb{}^J w$ is the unique decomposition with $w_J\in\We_J$, $\mb{}^J w\in\mb{}^J\We$.
\vspace*{1ex}\\
There are many important properties of a Coxeter system, its length function, its Bruhat order, and its minimal 
coset representatives. For all of this we also refer to the book \cite{Hu}. We only list four properties which we 
will use quite often:\vspace*{1ex}\\
{\bf 1)} The {\it Z-Lemma} of V. V. Deohdar, \cite{De}, Theorem 1.1 (II) : Let $c,d\in\We$. Let $s\in S$ such that $cs<c$ and $ds<d$. Then the following conditions are equivalent:\\
i) $c\leq d$. \\
ii) $cs\leq ds$.\\
iii) $cs\leq d$.\vspace*{1ex}\\
From the Z-Lemma follows easily:
Let $c,d\in\We$ such that $c\leq d$. Let $s\in S$. Then we have:\\
a) If $cs<c$ and $ds<d$, then $cs\leq ds$.\\
b) If $cs>c$ and $ds>d$, then $cs\leq ds$.\\
c) If $cs>c$ and $ds<d$, then $c\leq ds$ and $cs\leq d$.\vspace*{1ex}\\
Because we have either $ds<d$ or $ds>d\geq c$, and either $cs<c\leq d$ or $cs>c$, we get immediately from a) 
and b):\\
a') If $cs<c$, then $cs\leq ds$.\\
b') If $ds>d$, then $cs\leq ds$.\vspace*{1ex}\\
Furthermore we get by combining a) and c), and also b) and c):\\
d) If $ds<d$, then $cs\leq d$.\\
e) If $cs>c$, then $c\leq ds$.\vspace*{1ex}\\
{\bf 2)} The following Lemma of M. J. Dyer can be proved directly, or it can easily be obtained from the Z-Lemmma. 
In fact M. J. Dyer proved it directly for his more general twisted Bruhat orders to show the Z-Lemma for this 
orders. Compare the proof of \cite{Dy2}, Proposition 1.9: 
Let $t\in T$ and $s\in S$, $s\neq t$. Then:\\
a) If $z<zt$ then $zs<zts$.\\
b) If $zt< z$ then $zts< zs$.\vspace*{1ex}\\
{\bf 3)} The {\it subword property} of the Bruhat order of V. V. Deodhar, \cite{De}, Theorem 1.1 (III): 
Let $w\in\We\setminus\{1\}$ and fix a reduced expression $w=s_1s_2\cdots s_n$, $s_i\in S$. Then $v\leq w$ if and only 
if $v$ can be written as a subexpression of this
reduced expression of $w$, i.e., $v=1$ or $v=s_{i_1}s_{i_2}\cdots s_{i_k}$ where $1\leq i_1<i_2<\cdots <i_k\leq n$. Moreover it is possible to choose this subexpression reduced.\vspace*{1ex}\\
{\bf 4)} The alternatives given by V. V. Deohdar in \cite{De}, Lemma 3.1 and Lemma 3.2: 
Let $J\subseteq S$ and $w\in\We^J$. Let $s\in S$. Then exactly one of the following three cases holds:\\
a) $sw<w$. In this case $sw\in\We^J$.\\
b) $sw>w$ and $sw\in\We^J$.\\
c) $sw>w$ and $sw\notin \We^J$. In this case $sw=w\ti{s}$ for some $\ti{s}\in J$.
%
%
%
%
%
%
%
\section{Substitutes for multiplying $a\leq b$ by $w$\label{mult}}
If $a,b$ are elements of $\We$ such that $a\leq b$, and if $w$ is an arbitrary element of $\We$, it is not 
possible to conclude $aw\leq bw$. In this section we obtain to some extend substitutes for such a rule. These 
generalize Lemma 2.1 (i), (ii), and Lemma 2.2 of \cite{Pu3}, and show that these lemmas have a common origin.\\
The extended Bruhat orders, which will be introduced later, have many different characterizations. It is not 
obvious that these characterizations are equivalent. The results of this section are used in particular to show 
this equivalence. They are also used in this paper in the proofs of many other propositions and theorems.
\vspace*{1ex}\\ 
Let $a,b\in\We$. To cut short our notation we often denote the product $ab\in\We$ by
\begin{eqnarray*}
  \left\{\begin{array}{ccl}
  a\Box b & \mb{ if } & l(ab)=l(a)+l(b)\;, \\
  a\rhd b & \mb{ if } & l(ab)=l(a)-l(b)\;,\\
  a\lhd b & \mb{ if } & l(ab)=-l(a)+l(b)\;.
\end{array}\right.\;.
\end{eqnarray*}
These symbols are made in such a way that the length of an element with a vertical line as neighbour is 
counted positive. For $a,b,c\in\We$ we have
\begin{eqnarray*}\label{s1}
    abc\;=\; (a\Box b)\Box c  \quad\iff\quad  l(abc)=l(a)+l(b)+l(c)   
                                   \quad\iff\quad           abc\;=\; a \Box (b\Box c)\;,\\
   abc\; =\; (a\lhd b )\rhd c  \quad\iff\quad  l(abc)=-l(a)+l(b)-l(c)    
                                \quad\iff\quad            abc\;=\; a\lhd (b \rhd c)\;.\label{s2}
\end{eqnarray*}
In these cases we denote $abc$ by $a\Box b\Box c$, resp. $a\lhd b \rhd c $ for short. These properties are not difficult 
to check. For example suppose that $l(abc)=-l(a)+l(b)-l(c)$. Then $ab=a\lhd b$ because of
\begin{eqnarray*}
  -l(a)+l(b)\leq l(ab)=l(abcc^{-1})\leq l(abc)+l(c^{-1})=-l(a)+l(b)-l(c)+l(c^{-1})=-l(a)+l(b)\;.
\end{eqnarray*}
Now $abc = (a\lhd b )\rhd c$ follows from $l(abc)=-l(a)+l(b)-l(c)=l(a\lhd b)-l(c)$.\vspace*{1.5ex}\\
The next theorem is the main theorem of this section, describing to which extend it is possible to multiply 
$a\leq b$ by $w$. Part b) generalizes Lemma 2.1 (i) of \cite{Pu3}. It states that if $a,b,w\in\We$ 
such that $a\leq b$ and $aw=a\Box w$, then $a\Box w\leq b\Box w^+$ for some $w^+\leq w$.
\begin{Theorem}\label{P1} Let $a,b\in\We$ such that $a\leq b$. Let $w\in\We$.\vspace*{0.5ex}\\
{\bf a)} There exists an element $w^-\in\We$, $w^-\leq w$, such that 
$a\rhd w^- \leq bw$.\\
If in addition $aw=a\rhd w$ then $a\rhd w \leq bw$.\vspace*{0.5ex}\\
{\bf b)} There exists an element $w^+\in\We$, $w^+\leq w$, such that 
$aw\leq b\Box w^+$.\\
If in addition $bw=b\Box w$ then $aw\leq b\Box w$.
\end{Theorem}
\Proof The statement '$a\leq b$ and $bw=b\Box w$ implies $aw\leq b\Box w$' is an 
immediate consequence of the subword property of the Bruhat order. Nevertheless we derive 
it in another way, showing how it fits to the other statements of the theorem.\vspace*{1ex}\\
For the proof we use the conclusions of the Z-Lemma given in Section \ref{NP}. The letters a), b), c), a'), b'), d) 
and e) used in the proof refer to these conclusions.\vspace*{1ex}\\ 
If $w=1$ then $w^-:=1$ and $w^+:=1$ satisfy the required conditions. Now let $w\in\We\setminus\{1\}$ and fix a 
reduced expression $w=s_1s_2\cdots s_n$, $s_i\in S$.
Define recursively
\begin{eqnarray*}
  w_0^-  &:=& 1\;,\\
  w_k^-  &:=& \left\{\begin{array}{ccc}
   w_{k-1}^-s_k & \mb{if} & aw_{k-1}^-s_k< a w_{k-1}^-\\ 
   w_{k-1}^- & \mb{if} & aw_{k-1}^-s_k> a w_{k-1}^-   
\end{array}\right. ,\qquad k=1,2,\ldots,n\;,\\
\end{eqnarray*}
\begin{eqnarray*}
  w_0^+  &:=& 1\;,\\
  w_k^+  &:=& \left\{\begin{array}{ccc}
   w_{k-1}^+ & \mb{if} & bw_{k-1}^+s_k< b w_{k-1}^+\\ 
   w_{k-1}^+s_k & \mb{if} & bw_{k-1}^+s_k> b w_{k-1}^+   
\end{array}\right., \qquad k=1,2,\ldots,n\;.
\end{eqnarray*}
We show that $w^-:=w_n^-$ and $w^+:=w_n^+$ satisfy the required conditions:
$w_n^-$ and $w_n^+$ are obtained as subexpressions of the reduced expression $s_1s_2\cdots s_n$ of $w$. Therefore
$w_n^-\leq w$ and $w_n^+\leq w$.
Now set $w_0:=1$ and $w_k:=s_1s_2\cdots s_k$, $k=1,2,\ldots, n$. Inductively we show 
\begin{eqnarray*}
  a \rhd w_k^-\;\leq\; b w_k \quad \mb{ and }\quad a w_k\;\leq \;b \Box w_k^+  \quad,\quad k=0,1,\ldots,n\;.
\end{eqnarray*}
Since $a\leq b$ this is valid for $k=0$. Now suppose that it is valid for $k-1$, where $k\in\{1,\ldots,n\}$.\vspace*{1ex}\\
We have $aw_{k-1}^-\leq bw_{k-1}$. If $aw_{k-1}^-s_k< aw_{k-1}^-$ then due to a') we get 
$aw_{k-1}^-s_k \leq bw_{k-1}s_k$. If $aw_{k-1}^-s_k> aw_{k-1}^-$ then due to e) we get 
$ aw_{k-1}^-\leq bw_{k-1}s_k$. Due to our definition of $w_k^-$ and $w_k$ we have shown 
$aw_k^-\leq bw_k$.\vspace*{1ex}\\
We have $aw_{k-1}\leq bw_{k-1}^+$. If $bw_{k-1}^+s_k< bw_{k-1}^+$ then due to d) we get
$ aw_{k-1}s_k \leq bw_{k-1}^+ $.
If $bw_{k-1}^+s_k> bw_{k-1}^+$ then due to b') we get $aw_{k-1}s_k\;\leq\; bw_{k-1}^+s_k$.
Due to our definition of $w_k$ and $w_k^+$ we have shown $aw_k\leq bw_k^+$.\vspace*{1ex}\\
We have $l(aw_{k-1}^-)=l(a)-l(w_{k-1}^-)$. If $aw_{k-1}^-s_k< aw_{k-1}^-$ then
\begin{eqnarray*}
  l(a)-l(w_{k-1}^-)-1 \;=\; l(aw_{k-1}^-)-1\;=\; l(aw_{k-1}^-s_k)\;\geq\;l(a)-l(w_{k-1}^-s_k) \;\geq\;
  l(a)-l(w_{k-1}^-)-1 \;.
\end{eqnarray*}
Since $w_k^-=w_{k-1}^-s_k$ we get $l(aw_k^-)=l(a)-l(w_k^-)$.
If $aw_{k-1}^-s_k> aw_{k-1}^-$ then $w_k^-=w_{k-1}^-$. Trivially we get $l(aw_k^-)=l(a)-l(w_k^-)$.\vspace*{1ex}\\
We have $l(bw_{k-1}^+)=l(b)+l(w_{k-1}^-)$. If $bw_{k-1}^+s_k> bw_{k-1}^+$ then
\begin{eqnarray*}
  l(b)+l(w_{k-1}^+)+1 \;=\; l(bw_{k-1}^+)+1\;=\; l(bw_{k-1}^+s_k)\;\leq\;l(b)+l(w_{k-1}^+s_k) \;\leq\;
  l(b)+l(w_{k-1}^+)+1 \;.
\end{eqnarray*}
Since $w_k^+=w_{k-1}^+s_k$ we get $l(bw_k^+)=l(b)+l(w_k^+)$.
If $bw_{k-1}^+s_k< bw_{k-1}^+$ then $w_k^+=w_{k-1}^+$. Trivially we get $l(bw_k^+)=l(b)+l(w_k^+)$.\vspace*{1ex}\\
Now let $aw=a\rhd w$. Consider the elements
\begin{eqnarray*}
 a=aw_0,\quad aw_1,\quad aw_2,\quad\cdots,\quad aw_n=aw\;.
\end{eqnarray*}
At every step from $a w_{k-1}$ to $aw_{k}=aw_{k-1}s_k$, $k\in\{1,2,\ldots,n\}$, the length drops at most by one. To reach the value $l(aw)=l(a)-l(w)=l(a)-n$ the length has to drop at every step by one. Then from the definition of 
$w^-$ follows $w^-=w$.\vspace*{1ex}\\
Let $bw=b\Box w$. Consider the elements
\begin{eqnarray*}
 b=bw_0,\quad bw_1,\quad bw_2,\quad\cdots,\quad bw_n=bw\;.
\end{eqnarray*}
At every step from $b w_{k-1}$ to $bw_{k}=bw_{k-1}s_k$, $k\in\{1,2,\ldots,n\}$, the length increases at most by 
one. To reach the value $l(bw)=l(b)+l(w)=l(b)+n$ the length has to increase at every step by one. Then the 
definition of $w^+$ gives $w^+=w$.\\
\End
Lemma 2.2 of \cite{Pu3} states: If $I\subseteq S$, $a,b\in\We^I$, $u,v\in\We_I$ such that $au\leq bv$, then there exists elements $u_1,u_2\in\We$ such that $u=u_1\Box u_2$ and $au_1\leq b$ and $u_2\leq v$.
This Lemma is a particular case of part a) of the following Corollary of Theorem \ref{P1}, when 
applied to this situation. (Use also $l(w)=l(w^I)+l(w_I)$ for $w\in\We$.)
\begin{Cor}\label{P2} Let $a,b\in\We$ and $v\in\We$. Then:\vspace*{0.5ex}\\
{\bf a)} If $a \leq bv $ then there exists an element $\ti{v}\in\We$, $\ti{v}\leq v$ such that 
$a\rhd \ti{v}^{-1}\leq b $.\\
If we have even $a \leq b\rhd v $ then $a v^{-1}\leq b $.\vspace*{0.5ex}\\
{\bf b)} If $av\leq b$ then there exists an element $\ti{v}\in\We$, $\ti{v}\leq v$ such that 
$a\leq b\Box \ti{v}^{-1}$.\\
If we have even $a\Box v \leq b$ then $a\leq b v^{-1}$.
\end{Cor}
\Proof
Applying part a) of Theorem \ref{P1} to the inequality $a\leq bv$, setting $w=v^{-1}$, we find an element $\ti{v}\leq v$ 
such that $a\rhd\ti{v}^{-1}\leq bv v^{-1}=b$.
If $bv=b\rhd v$ then $b=(b\rhd v)\Box v^{-1}$ because of
\begin{eqnarray*}
   l(b)\;=\; l((b\rhd v)v^{-1})\;\leq\;l(b\rhd v)+l(v^{-1})\;=\;l(b)-l(v)+l(v^{-1})\;=\;l(b)\;.
\end{eqnarray*}
Applying part b) 'In addition ...' of Theorem \ref{P1} to the inequality $a\leq bv$, setting $w=v^{-1}$, we find 
$av^{-1}\leq b$.\vspace*{1ex}\\
Applying part b) of Theorem \ref{P1} to the inequality $av\leq b$, setting $w=v^{-1}$, we find an element 
$\ti{v}\leq v$ such that $a= a v v^{-1}\leq b\Box \ti{v}^{-1}$. 
If $av=a\Box v$ then $a=(a\Box v)\rhd v^{-1}$. Applying part a) 'In addition ...' of Theorem \ref{P1} to the 
inequality $av\leq b$, setting $w=v^{-1}$, we find $a\leq bv^{-1}$.\\
\End
The following canceling rule follows immediately from the statements 'If we have even ...' in Corollary \ref{P2}. 
Part b) generalizes Lemma 2.1 (ii) of \cite{Pu3}. It states that for $a,b,w\in\We$ the inequality 
$a\Box w\leq b\Box w$ implies $a\leq b$.
\begin{Cor}\label{P3} Let $a,b\in\We$ and $w\in\We$. Then:\vspace*{0.5ex}\\
{\bf a)} If $aw\leq b \rhd w$ then $a\leq b$.\vspace*{0.5ex}\\
{\bf b)} If $a\Box w\leq bw $ then $a\leq b$.
\end{Cor}
The following conclusion of Corollary \ref{P2} will be used in the next section several times to show the equivalence 
of different characterizations of the extended Bruhat orders.
\begin{Cor}\label{P4} Let $a,b\in\We$. Let $S(\;\; )$ be a statement about the elements $\We$, such that if 
$S(w)$ is true for an element $w\in\We$, then also $S(\ti{w})$ is true for all elements $\ti{w}\in\We$, 
$\ti{w}\leq w$. Then the following four statements are equivalent:\vspace*{0.5ex}\\
(i) There exists an element $v\in\We$ such that $a\leq bv$ and $S(v)$.\vspace*{0.5ex}\\
(ii) There exists an element $v\in\We$ such that $av^{-1}\leq b$ and $S(v)$.\vspace*{0.5ex}\\
(i') There exists an element $v\in\We$ such that $a\leq b\Box v$ and $S(v)$.\vspace*{0.5ex}\\
(ii') There exists an element $v\in\We$ such that $a\rhd v^{-1}\leq b$ and $S(v)$.
\end{Cor}
\Proof Obviously (i') implies (i), and (ii') implies (ii). Due to part a) of Corollary \ref{P2}, (i) implies (ii'). 
Due to part b) of Corollary \ref{P2}, (ii) implies (i').\\
\End
The theorems and corollaries in this section are substitutes for multiplying $a\leq b$ by $w$ from the right. 
Clearly there are also the corresponding versions for multiplying $a\leq b$ by $w$ from the left. These can be 
obtained by applying the inverse map. If we quote a theorem or corollary of this section in this paper, it means 
we refer to both versions.  
%
%
%
%
%
%
%
%
%
%
%
%
\section{The extended Bruhat orders and length functions\label{BL}}
For the rest of the paper we fix a subset $N$ of $S$, and a component $C$ of $N$, (i.e., 
$C\subseteq N$ and $s\ti{s}=\ti{s}s$ for all $s\in N\setminus C$ and $\ti{s}\in C$).\vspace*{1ex}\\
In this section we introduce the extended Bruhat orders $\leq_{\eps\delta}$, and compatible extended length 
functions $l_{\eps\delta}$, $\eps,\delta\in\{+,-\}$, on a certain set $\We(N,C)$.
We give several characterizations of the extended Bruhat orders. We investigate if there exist isomorphisms or 
anti-isomorphisms.\vspace*{1ex}\\ 
We denote by $\We^{op}$ the opposite group of $\We$. Equip the group $\We\times\We^{op}$ with the action on itself 
by left multiplication, i.e.,
\begin{eqnarray*}
   (u,v)(a,b):= (ua,bv) \quad\mb{ where }\quad u,v,a,b \in\We.
\end{eqnarray*} 
Equip the group $\We\times\We^{op}$ with the involution of groups $\mb{}^{inv}:\We\times\We^{op} \to\We\times\We^{op}$ 
given by
\begin{eqnarray*}
 (a,b)^{inv}:=(b^{-1}, a^{-1}) \quad\mb{ where }\quad a,b\in\We.
\end{eqnarray*}
The subgroups
\begin{eqnarray*}
  \We_C\times 1 \quad,\quad \Mklz{(v,v^{-1})}{v\in\We_{N\setminus C}}\quad,\quad 1\times (\We_C)^{op}
\end{eqnarray*}
of $\We\times\We^{op}$ commute pairwise, because every simple reflections of $C$ commutes with every simple 
reflection of $N\setminus C$. Therefore 
the product of these three subgroups is again a subgroup. It is also invariant under the involution $\mb{}^{inv}$. 
\begin{Def}\label{BL1} Let $\We(N,C)$ be the quotient of the group $\We\times\We^{op}$ by the subgroup 
\begin{eqnarray*}
   \bigl(\We_C\times 1\,\bigr)\,\Mklz{(v,v^{-1})}{v\in\We_{N\setminus C}}\,\bigl(\,1\times (\We_C)^{op}\bigr)\;.
\end{eqnarray*}
Denote the image of $(a,b)\in\We\times\We^{op}$ under the canonical projection by $a\idem b$.\vspace*{1ex}\\
Equip $\We(N,C)$ with the descended $\We\times\We^{op}$-action, i.e.,
\begin{eqnarray*}
   (u,v)\, a\idem b:= ua\idem bv \quad\mb{ where }\quad u,v,a,b \in\We.
\end{eqnarray*}
Equip $\We(N,C)$ with the descended involution, also denoted by $\mb{}^{inv}$, i.e.,
\begin{eqnarray*}
 (a\idem b)^{inv}:= b^{-1}\idem a^{-1} \quad\mb{ where }\quad a,b\in\We.
\end{eqnarray*}
\end{Def}
{\bf Remark:} The $\We\times\We^{op}$-module $\We(S,\emptyset)$ identifies with the 
$\We\times\We^{op}$-module $\We$, the involution $\mb{}^{inv}$ of $\We(S,\emptyset)$ with the inverse map of 
$\We$.\vspace*{1ex}\\
We call $N(\idem) :=\Mklz{w\in\We}{w\idem =\idem w}$ the normalisator, we call 
$C_L(\idem) := \Mklz{w\in\We}{w\idem =\idem  }$ and $C_R(\idem) :=  \Mklz{w\in\We}{\idem w=\idem  }$ the 
left and right centralizators of $\idem\in\We(N,C)$. 
\begin{Prop}\label{BL2} We have $N(\idem)=\We_N$, $C_L(\idem) =C_R(\idem) =\We_C$.
\end{Prop}
\Proof As an example we show $N(\idem)=\We_N$. The statements about the left and right centralizators are 
shown in a similar way. Let $w\in\We$. By definition $w\idem =\idem w$ if and only if there exist elements
$u\in\We_C$, $v\in\We_{N\setminus C}$, and $\ti{u}\in\We_C$ such that 
\begin{eqnarray*}
  w=1uv \quad \mb{ and }\quad 1=v^{-1}\ti{u}w\;.
\end{eqnarray*}
This is equivalent to $w\in\We_C\We_{N\setminus C}=\We_N$.\\
\End
The elements of $\We(N,C)$ can be represented in particular ways, which will be very useful:
\begin{Prop}\label{BL3} Let $x\in\We(N,C)$. Then:\vspace*{0.5ex}\\
(I) There exist uniquely determined elements $a\in\We^C$, $b\in \mb{}^N\We$ such that 
\begin{eqnarray}\label{enf1}
    x\;=\;a\idem b \;.
\end{eqnarray}
(II) There exist uniquely determined elements $a\in\We^N$, $b\in \mb{}^C \We$ such that
\begin{eqnarray}\label{enf2}
    x\;=\;a\idem b \;.
\end{eqnarray}
(III) There exist uniquely determined elements $a\in\We^N$, $c\in\We_{N\setminus C}$, $b\in \mb{}^N \We$ such that 
\begin{eqnarray}\label{enf3}
    x\;=\;a c\idem b \;=\;a\idem c b\;.
\end{eqnarray}
We call the expression in (\ref{enf1}), (\ref{enf2}), resp. (\ref{enf3}) on the right the normal 
forms I, II, resp. III of $x$. 
\end{Prop}
{\bf Remarks:} 
{\bf (1)} By applying the involution $\mb{}^{inv}:\We(N,C)\to\We(N,C)$ to an element $x$ in normal form I resp. II 
we obtain the element $x^{inv}$ in 
normal form II resp. I. By applying this map to an element $x$ in normal form III we obtain the element 
$x^{inv}$ in normal form III.\\ 
{\bf (2)} Because the multiplication map of $\We$ restricts to bijective maps  
\begin{eqnarray*}
    \We^N\times \We_{N\setminus C}\;\to\; \We^C \qquad \mb{and}\qquad 
    \We_{N\setminus C}\times \mb{}^N\We\;\to\;\mb{}^C\We\;\;,
\end{eqnarray*} 
we can immediately read off the normal forms I and II from the normal form III.\vspace*{1ex}\\ 
\Proof Due to these remarks it is sufficient to show (I). To show the existence of normal form I let 
$c\idem d\in\We(N,C)$, $c,d\in\We$. By using Proposition \ref{BL2} we get
\begin{eqnarray*}
   c\idem d\;=\; c \idem d_{N\setminus C}d_C \,\mb{}^N d\;=\; c d_{N\setminus C} \idem  \,\mb{}^N d\;=\;
  (c d_{N\setminus C})^C \idem \, \mb{}^N d\;.
\end{eqnarray*}
To show the uniqueness of normal form I let $a\idem b= \ti{a}\idem \ti{b}$ with $a,\ti{a}\in \We^C$ and 
$b,\ti{b}\in \mb{}^N\We$. Then by Definition \ref{BL1} there exist elements $u\in\We_C$,  
$v\in\We_{N\setminus C}$, and $w\in\We_C$ such that
\begin{eqnarray*}
  \ti{a}= a uv \qquad\mb{ and }\qquad \ti{b}=v^{-1}w b\;.
\end{eqnarray*}
Since $b,\ti{b}\in \mb{}^N\We$ and $v^{-1}w\in\We_N$ the second equation implies $\ti{b}= b$ and $v^{-1}w=1$. 
Since $v\in\We_{N\setminus C}$ and $w\in\We_C$ it follows $v=1$. Inserting in the first 
equation we get $\ti{a}= a u$. Since $a,\ti{a}\in \We^C$ and $u\in\We_C$ we find $\ti{a}= a$.\\
\End
Let $c,d\in\We^N$ and $w\in\We_N$ such that $cw\leq d$. Then $cw=c\Box w$ and from Theorem \ref{P1} 'In addition ... ' we find $c\ti{w}\leq c w\leq d$ for all $\ti{w}\in\We$, $\ti{w}\leq w$. Similar things hold if $c w^{-1}\leq d$, if $wc\leq d$, 
or if $w^{-1}c\leq d$. Therefore by Corollary \ref{P4} the relations given in the following definition are well 
defined.
\begin{Def}\label{BL4} Let $\eps,\delta\in\{+,-\}$. Define a relation $\leq_{\eps\delta}$ on $\We(N,C)$, which we 
call an extended Bruhat order, as follows:
For $x_1, x_2\in\We(N,C)$ let $x_1=a_1 \idem c_1 b_1= a_1 c_1\idem b_1$, $x_2=a_2 \idem c_2 b_2= a_2 c_2\idem b_2$ be its normal forms III. Set $x_1\leq_{\eps\delta}x_2$ if there exist 
$u,v\in\We_{N\setminus C}$ such that
\begin{eqnarray*}
  \left\{ \begin{array}{ccl}
    a_1 u^{-1}\leq a_2    & \mb{for} & \delta=+\\
    a_1\geq a_2 u     & \mb{for} & \delta=-
  \end{array}\right\} 
\quad\mb{ and }\quad (x) \quad \mb{ and }\quad
\left\{ \begin{array}{ccl}
    b_1\geq v b_2    & \mb{for} & \eps=+\\
    v^{-1} b_1 \leq b_2     & \mb{for} & \eps=-
  \end{array}\right\} 
\end{eqnarray*}
holds. Here for $(x)$ any of the following eight statements can be taken:
\begin{eqnarray*}
\begin{array}{ccc}
(i)   & &  c_1\leq u^{-1}c_2 v^{-1}\\
(ii)  & &    u c_1 \leq c_2 v^{-1}  \\
(iii) & &  u c_1 v\leq c_2 \\
(iv)  & &  c_1 v \leq u^{-1}c_2  \vspace*{0.5ex}\\ 
(i')  & &  c_1 \leq  u^{-1}\Box c_2 \Box v^{-1} \\
(ii') & &  u\lhd c_1  \leq  c_2\Box v^{-1} \\
(iii')& &  u \lhd c_1 \rhd v  \leq  c_2\\
(iv') & &  c_1\rhd v \leq  u^{-1}\Box c_2 
\end{array}
\end{eqnarray*}
\end{Def}
{\bf Remark:} Identify $\We(N,C)$ with $\We^N\times \We_{N\setminus C}\times \mb{}^N\We$ as a set. Then the 
restriction of the relation $\leq_{\eps\delta}$ to one of this factors is always the Bruhat order or the 
inverse Bruhat order:
\begin{eqnarray*}
  \begin{array}{l|c|c|c|c|}
      &\quad\qquad -+\qquad\quad &\quad\qquad ++\qquad\quad & \quad\qquad--\qquad\quad &\quad\qquad +-\qquad\quad \\
     \hline 
     \mb{first factor}  &\leq  & \leq                 & \mb{inverse of }\leq & \mb{inverse of }\leq \\
     \mb{middle factor} &\leq  & \leq                 & \leq                 & \leq \\
     \mb{last factor}   &\leq  & \mb{inverse of }\leq & \leq                 & \mb{inverse of }\leq 
  \end{array}
\end{eqnarray*}
For $\We(S,\emptyset)$ identified with $\We$ the four extended Bruhat orders coincide with the Bruhat order 
on $\We$.
\begin{Prop}\label{BL5} For $\eps,\delta\in\{+,-\}$ the relation $\leq_{\eps\delta}$ is a partial order on $\We(N,C)$.
\end{Prop}
\Proof We only show this for $\leq_{++}$. The proofs for the other three relations are similar.
Obviously the relation $\leq_{++}$ is reflexive. To show that it is anti-symmetric let 
$x_1,x_2\in\We(N,C)$ such that $x_1\leq_{++} x_2$ and $x_2\leq_{++}x_1$. Let $x_1=a_1 c_1 e b_1$ and 
$x_2=a_2 c_2 e b_2$ be the normal forms III. By definition of the relation $\leq_{++}$ there exist elements 
$u,\ti{u},v,\ti{v}\in\We$ such that
\begin{eqnarray}
     a_1 u^{-1} \leq a_2\;,   &&  a_2 \ti{u}^{-1} \leq a_1 \;,           \label{as1}\\
     c_1\leq u^{-1}c_2 v^{-1}\;, && c_2\leq \ti{u}^{-1}c_1\ti{v}^{-1}\;, \label{as2}\\
     b_1\geq v b_2 \;,   && b_2\geq \ti{v}b_1\;.                         \label{as3}
\end{eqnarray}
Since $a_1 u^{-1}=a_1\Box u^{-1}$ we have $a_1\leq a_1 u^{-1}$. Similarly $a_2\leq a_2 \ti{u}^{-1}$. Together with 
(\ref{as1}) it follows $a_1=a_2$ and $u=\ti{u}=1$. In the same way from (\ref{as3}) we get $b_1=b_2$ and 
$v=\ti{v}=1$. Inserting $u=\ti{u}=v=\ti{v}=1$ in (\ref{as2}) we find $c_1=c_2$.\vspace*{1ex}\\  
To show the transitivity let $x_1,x_2,x_3\in\We(N,C)$ such that $x_1\leq_{++} x_2$ and $x_2\leq_{++}x_3$. Let $x_1=a_1 c_1 e b_1$, $x_2=a_2 c_2 e b_2$, and $x_3=a_3 c_3 e b_3$ be the normal forms III. By definition of the relation 
$\leq_{++}$ there exist elements $u,\ti{u},v,\ti{v}\in\We$ such that
\begin{eqnarray}
     a_1 u^{-1} \leq a_2\;,   &&  a_2 \ti{u}^{-1} \leq a_3 \;,           \label{t1}\\
     c_1\leq u^{-1}c_2 v^{-1}\;, && c_2\leq\ti{u}^{-1}c_3\ti{v}^{-1}\;, \label{t2}\\
     b_1\geq v b_2 \;,   && b_2\geq \ti{v}b_3\;.                         \label{t3}
\end{eqnarray}
Applying two times Theorem \ref{P1} b) to $c_2\leq\ti{u}^{-1}c_3\ti{v}^{-1}$ we find elements $u^+,v^+\in\We$, 
$u^+\leq u$ and $v^+\leq v$ such that $u^{-1}c_2 v^{-1}\leq(u^+)^{-1}\Box \ti{u}^{-1}c_3\ti{v}^{-1} 
\Box (v^+)^{-1}$. By the first inequality of (\ref{t2}) it follows $c_1\leq (\ti{u}u^+)^{-1}c_3(v^+\ti{v})^{-1}$.
Since $a_1u^{-1}=a_1\Box u^{-1}$ and $(u^+)^{-1}\leq u^{-1}$ we get from Theorem \ref{P1} b) 
$a_1 (u^+)^{-1} \leq a_1 u^{-1}$. With the first inequality of (\ref{t1}) it follows 
$a_1 (u^+)^{-1} \leq a_2$. Now $a_2 \ti{u}^{-1}=a_2 \Box\ti{u}^{-1}$. Applying Theorem \ref{P1} b) once more we 
get $a_1(u^+)^{-1}\ti{u}^{-1}\leq a_2\ti{u}^{-1}$. With the second inequality of (\ref{t1}) it follows 
$a_1(\ti{u}u^+)^{-1}\leq a_3$.
Similarly from (\ref{t3}) we get $b_1\geq v^+\ti{v}b_3$. By the definition of the relation $\leq_{++}$ we have 
shown $x_1\leq_{++}x_3$.\\
\End
\begin{Def}\label{BL6} Let $\eps,\delta\in\{+,-\}$. Define a function $l_{\eps\delta}:\We(N,C)\to\Z$,
which we call an extended length function, as follows: For $x\in\We(N,C)$ let $x=a\idem c b= a c\idem b$ be its 
normal form III. Set 
\begin{eqnarray*}
 l_{\eps\delta}(x) &:=& \delta\, l(a) + l(c)-\eps\,l(b)\;.
\end{eqnarray*}
\end{Def}
{\bf Remark:} Identify $\We(N,C)$ with $\We^N\times \We_{N\setminus C}\times \mb{}^N\We$ as a set. Then the 
restriction of the extended length function $l_{\eps\delta}$ to one of these factors is always the length function or 
the negative of the length function. It matches with the restrictions of the extended Bruhat order $\leq_{\eps\delta}$:
\begin{eqnarray*}
  \begin{array}{l|c|c|c|c|}
      &\quad\qquad -+\qquad\quad &\quad\qquad ++\qquad\quad & \quad\qquad--\qquad\quad &\quad\qquad +-\qquad\quad \\
     \hline 
     \mb{first factor}  & l &  l & -l & -l \\
     \mb{middle factor} & l &  l &  l &  l \\
     \mb{last factor}   & l & -l &  l & -l
  \end{array}
\end{eqnarray*}
For $\We(S,\emptyset)$ identified with $\We$ the four extended length functions coincide with the length function 
on $\We$.\vspace*{1ex}\\ 
Equip $\Z$ with its natural order. The next proposition shows that the extended length functions are compatible with 
the extended Bruhat orders.
\begin{Prop}\label{BL7} Let $\eps,\delta\in \{+,-\}$. Let $x,y\in\We(N,C)$ such that $x<_{\eps\delta} y$. Then 
also $l_{\eps\delta}(x)<l_{\eps\delta}(y)$. 
\end{Prop}
\Proof We only show this for $\leq_{++}$ and $l_{++}$. The proofs for the other three extended Bruhat orders 
and length functions are similar.
Let $x=a_1 c_1 e b_1$ and $y=a_2 c_2 e b_2$ be the normal forms III of $x,y$. By definition of the relation 
$\leq_{++}$ there exist elements $u,v\in\We$ such that
\begin{eqnarray}\label{lineq}
     a_1 u^{-1} \leq a_2  \quad,\quad   c_1\leq u^{-1}\Box c_2 \Box v^{-1}  \quad,\quad  
     b_1\geq v b_2 \;.           
\end{eqnarray}
Suppose that none of these three inequalities is proper. Since $a_1, a_2\in\We^N$ and $u\in\We_{N\setminus C}$ 
from the equation $a_1 u^{-1}=a_2 $ follows $a_1=a_2$ and $u=1$. In the same way from $ b_1=v b_2$ follows 
$b_1=b_2$ and $v=1$. Inserting in $c_1=u^{-1} c_2 v^{-1}$ we get $c_1=c_2$. Therefore we would have $x=y$ which is 
not possible.
Now by the length inequalities corresponding to the inequalities (\ref{lineq}), and by
$l(a_1u^{-1})=l(a_1)+l(u^{-1})$ and $l(v b_2)=l(v)+l(b_2)$ we find
\begin{eqnarray*}
   l_{++}(x) \;=\; l(a_1)+ l(c_1) -l(b_1)\;=\;l(a_1\Box u^{-1}) -l(u) +l(c_1) -l(b_1)\;< \;\\
   l(a_2) -l(u) +l(u^{-1}\Box c_2 \Box v^{-1}) -l(v\Box b_2)\;=\;l(a_2)+l(c_2) -l(b_2)
   \;=\;l_{++}(y)\;.
\end{eqnarray*}
\End 
To complete the elementary properties of the extended Bruhat orders we investigate next which of them are isomorphic 
or anti-isomorphic. After that we investigate which of the extended Bruhat orders can also be characterized by using 
normal form I or II.\\
To analyze the set $\We(e)=\We \idem\We$ equipped with the Bruhat-Chevalley order for a finite Coxeter group $\We$, 
M. S. Putcha defines in \cite{Pu3} two posets $\We_{N,C}$ and $\We_{N,C}^*$. 
In Theorem 2.5 i), iii) of \cite{Pu3} he shows that $\We(e)$ and $\We_{N,C}$ are isomorphic, and 
$\We_{N,C}$ and $\We_{N,C}^*$ are anti-isomorphic.\\ 
The poset $\We(e)$ identifies with ($\We(N,C)$,$ \leq_{++}$), but $\leq_{++}$ defined differently by using normal 
form I. 
The definition of the poset $\We_{N,C}$ uses the longest element of ($\We_{N\setminus C}$, $N\setminus C$). 
It is none of the descriptions given in this paper. 
The poset $\We_{N,C}^*$ identifies with ($\We(N,C)$, $\leq_{-+}$).
The concatenation $\We(e)\to\We_{N,C}\to\We_{N,C}^*$ of the isomorphism and anti-isomorphism given in the proofs is the 
left multiplication by the longest element of ($\We$, $S$).\\ 
The main parts of the proofs of the following Theorem \ref{BL9} and Proposition \ref{BL10} are obtained from the proof of 
Theorem 2.5 i), iii) of \cite{Pu3} by eliminating the intermediate $\We_{N,C}$, and by isolating and removing 
the transformation of the normal forms.\vspace*{1ex}\\  
The following compatibilities between the involution $\mb{}^{inv}:\We(N,C)\to\We(N,C)$ and the extended Bruhat orders 
and length functions follow immediately from its definitions and from Remark (1) after Proposition \ref{BL3}.
\begin{Prop}\label{BL8} Consider the involution $\mb{}^{inv}:\We(N,C)\to\We(N,C)$.\mb{}\\
{\bf a)} It is an isomorphism of $(\We(N,C),\leq_{++}, l_{++})$ and $(\We(N,C),\leq_{--}, l_{--})$.\\
{\bf b)} It is an automorphism of $(\We(N,C),\leq_{-+}, l_{-+})$.\\
{\bf c)} It is an automorphism of $(\We(N,C),\leq_{+-}, l_{+-})$.
\end{Prop}
In general each two of the extended Bruhat orders $\leq_{++}$, $\leq_{-+}$, and $\leq_{+-}$ are not isomorphic and 
not anti-isomorphic. This can be seen by looking at the smallest and biggest elements of these orders in case of 
a Coxeter system $(\We,S)$ with a subset $N\subseteq S$ such that $\We^N$ and $\We_N$ are  infinite:\\
$\bullet$ There is no smallest element and no biggest element of $\leq_{++}$.\\
$\bullet$ The smallest element of $\leq_{-+}$ is $\idem$, but there is no biggest element.\\
$\bullet$ If $\We_{N\setminus C}$ is infinite, then there is no smallest element and no biggest element 
of $\leq_{+-}$.\\
Let $\We_{N\setminus C}$ be finite and $u_0$ be the longest element of ($\We_{N\setminus C}$, $N\setminus C$). Then 
there is no smallest element of $\leq_{+-}$, but $u_0\idem $ is the biggest element.\vspace*{1ex}\\
The situation is different for a finite Coxeter group:
\begin{Theorem}\label{BL9} Let $\We$ be finite and $w_0$ be the longest element of $(\We, S)$. Denote by\\ 
\mb{} $\quad\Phi_{w_0 1}:\We(N,C)\to\We(N,C)$ the left multiplication by $w_0$,\\ 
\mb{} $\quad\Phi_{1 w_0}:\We(N,C)\to\We(N,C)$ the right multiplication by $w_0$, \\
\mb{} $\quad\Phi_{w_0 w_0}:\We(N,C)\to\We(N,C)$ the both sided multiplication by $w_0$. \\
As indicated by the lines in the following diagrams the maps $\Phi_{w_0 1}$, $\Phi_{1 w_0 }$ are involutive 
anti-isomorphisms, the map $\Phi_{w_0 w_0}$ is an involutive isomorphism between certain extended Bruhat orders:\\\\
\unitlength1.7ex
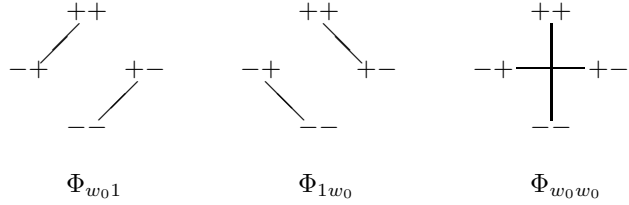
\hspace*{18ex}\begin{picture}(35,9)
\put(3,3){$--$}\put(3,9){$++$}\put(0,6){$-+$}\put(6,6){$+-$}\put(3,0){$\Phi_{w_0 1}$} 
\put(4.7,3.7){\line(1,1){2}}\put(1.7,6.7){\line(1,1){2}}
\put(15,3){$--$}\put(15,9){$++$}\put(12,6){$-+$}\put(18,6){$+-$}\put(15,0){$\Phi_{1 w_0 }$} 
\put(15.3,3.7){\line(-1,1){2}}\put(18.3,6.7){\line(-1,1){2}}
\put(27,3){$--$}\put(27,9){$++$}\put(24,6){$-+$}\put(30,6){$+-$}\put(27,0){$\Phi_{w_0 w_0 }$} 
\put(28.1,3.7){\line(0,1){5}}\put(26.3,6.4){\line(1,0){3.5}}
\end{picture}
\end{Theorem}
\Proof It is easy to check that $\Phi_{1 w_0}=\mb{}^{inv} \circ\Phi_{w_0 1}\circ \mb{}^{inv}$ and 
$\Phi_{w_0 w_0}=\Phi_{w_0 1}\circ \Phi_{1 w_0}=\Phi_{1 w_0}\circ \Phi_{w_0 1} $. Using these relations, the 
statements for $\Phi_{1 w_0}$ and $\Phi_{w_0 w_0}$ follow from the statements for $\Phi_{w_0 1}$.\\
As an example we show that $\Phi_{w_0 1}$ is an involutive anti-isomorphism from $(\We(N,C),\leq_{++})$ to 
$(\We(N,C),\leq_{-+})$, the proof of the remaining statement is similar.\vspace*{1ex}\\
Clearly the map $\Phi_{w_0 1}$ is involutive, in particular it is bijective.
Now let $v_0$ be the longest element of $\We_N$. Let $u_0$ be the longest element of $\We_{N\setminus C}$.
Let $t_0$ be the longest element of $\We_C$. Then $v_0=u_0 t_0=t_0 u_0$, and if $w\in\We^N$ then also 
$w_0wv_0\in\We^N$.
If $ac\idem b$ is the normal form III of an element of $\We(N,C)$, then 
\begin{eqnarray*}
    \Phi_{w_0 1}(ac\idem b)\;=\;w_0 a c \idem b\;=\;w_0 a (v_0  u_0 t_0) c \idem b\;=\;
    (w_0 a v_0)(u_0 c)(t_0\idem )b\;=\;(w_0 a v_0) (u_0 c) \idem b\;.
\end{eqnarray*}
Here the last expression is the normal form III of $\Phi_{w_0 1}(ac\idem b)$.\vspace*{1ex}\\ 
Let $a_1 c_1\idem b_1$, $a_2 c_2\idem b_2$ be the normal forms III of two elements of $\We(N,C)$. By definition  
$a_1 c_1\idem b_1 \leq_{++} a_2 c_2\idem b_2$ if and only if there exist elements $u,v\in\We_{N\setminus C}$ such 
that
\begin{eqnarray}\label{ineqw_0}
    a_1 u^{-1}\leq a_2  \quad\mb{ and }\quad c_1\leq u^{-1}c_2 v^{-1} \quad \mb{ and }\quad  b_1\geq v b_2\;\;.
\end{eqnarray}
We have $a_2 u v_0= a_2\Box (u v_0)$. Due to Theorem \ref{P1} b) the first inequality of (\ref{ineqw_0}) implies
\begin{eqnarray}\label{a1v0}
   a_1 v_0\;=\; (a_1 u^{-1}) u v_0 \;\leq\; a_2 u v_0\;\;.
\end{eqnarray}
Since $(a_1 u^{-1})(u v_0)=(a_1 u^{-1})\Box (u v_0)$ we can get back the the first inequality of (\ref{ineqw_0}) from (\ref{a1v0}) by applying the canceling rule Corollary \ref{P3} b).
Multiplying by $w_0$ from the left, reversing the order, inequality (\ref{a1v0}) is equivalent to
\begin{eqnarray}\label{ineqw_2}
   w_0 a_1 v_0 \;\geq\; w_0 a_2 u v_0\;=\;w_0 a_2 v_0 (v_0 u v_0)\;=\;(w_0 a_2 v_0) (u_0 u u_0) \;.
\end{eqnarray}
Multiplying by $u_0$ from the left, reversing the order, the second inequality of (\ref{ineqw_0}) is equivalent to
\begin{eqnarray}\label{ineqw_1}
   u_0 c_1\;\geq\; u_0u^{-1}c_2 v^{-1}\;=\;  (u_0 u u_0)^{-1} (u_0 c_2) v^{-1}\;.
\end{eqnarray}
By (\ref{ineqw_2}), (\ref{ineqw_1}), and the third inequality of (\ref{ineqw_0}) we have shown 
$a_1 c_1\idem b_1 \leq_{++} a_2 c_2\idem b_2$ if and only if 
$\Phi_{w_0 1}(a_2c_2\idem b_2)\leq_{-+}\Phi_{w_0 1}(a_1c_1\idem b_1)$.\\
\End
The extended Bruhat orders have been defined by using normal form III. As the following propositions show, 
there is also the possibility to characterize $\leq_{-+}$, $\leq_{++}$ in an easy way by using normal form I, 
and $\leq_{-+}$, $\leq_{--}$ by using normal form II. This is not possible 
for $\leq_{+-}$. It has the following reason: Take for example normal form I. If 
$\ti{a}_1 c_1 \idem b_1=\ti{a}_1 \idem c_1 b_1 $ is the normal form III of an element of $\We(N,C)$, we get normal 
form I by multiplying together the first and middle factor, i. e., $(\ti{a}_1 c_1) \idem b_1$.
Now identify $\We(N,C)$ with $\We^N\times \We_{N\setminus C}\times \mb{}^N\We$ as a set. To be able to transform 
the definition of an extended Bruhat order in a characterization with normal form I, the restrictions of the extended 
Bruhat order to the first and middle factor have to be uniform, i.e., on both terms the Bruhat order. This is only the 
case for $\leq_{-+}$, $\leq_{++}$.
\begin{Prop}\label{BL10} Let $\eps\in\{+,-\}$. Let $x_1, x_2\in\We(N,C)$ and let $x_1=a_1 \idem b_1$, 
$x_2=a_2 \idem b_2$ be its normal forms I. Then $x_1\leq_{\eps+}x_2$ if and only if there exists an element 
$v\in\We_{N\setminus C}$ such that
\begin{eqnarray*}
  (x) \quad\mb{ and }\quad  \left\{\begin{array}{ccl} 
                           b_1\geq v b_2 & for & \eps=+\\
                           v^{-1} b_1\leq b_2    & for & \eps=-
\end{array}\right\}
\end{eqnarray*}
holds. Here for $(x)$ any of the following four statements can be taken:
\begin{eqnarray*}
\begin{array}{lcc}
(i)  & & a_1\leq a_2 v^{-1} \\
(ii) & & a_1 v \leq a_2 \vspace*{0.5ex}\\
(i') & &  a_1\leq a_2 \Box v^{-1} \\
(ii') & & a_1 \rhd v \leq a_2 
\end{array}
\end{eqnarray*}
\end{Prop}
\Proof By Corollary \ref{P4} the different characterizations stated in this proposition are equivalent. Let 
$x_1=\ti{a}_1c_1\idem b_1$ and $x_2=\ti{a}_2c_2\idem b_2$ be the normal forms III of $x_1$ and $x_2$.
By definition $x_1\leq_{\eps +} x_2$ if there exist elements $u,v\in\We_{N\setminus C}$ such that
\begin{eqnarray*}
   \ti{a}_1 u^{-1}\leq \ti{a}_2 \quad\mb{ and }\quad c_1\leq u^{-1}c_2 v^{-1} \quad\mb{ and }\quad
   \left\{\begin{array}{ccl} 
                           b_1\geq v b_2 & for & \eps=+\\
                           v^{-1} b_1\leq b_2    & for & \eps=-
\end{array}\right\}\;.
\end{eqnarray*}
Applying Theorem \ref{P1} b) to $c_1\leq u^{-1}c_2 v^{-1}$ and $ \ti{a}_1 (u^{-1}c_2 v^{-1})=\ti{a}_1\Box u^{-1}c_2 v^{-1}$ we get $\ti{a}_1 c_1\leq \ti{a}_1u^{-1}c_2 v^{-1}$. Applying Theorem \ref{P1} b) to $\ti{a}_1 u^{-1}\leq \ti{a}_2 $ 
and $\ti{a}_2 (c_2 v^{-1})=\ti{a}_2\Box c_2 v^{-1}$  we get 
$\ti{a}_1u^{-1}c_2 v^{-1}\leq\ti{a}_2 c_2 v^{-1}$. It follows $ \ti{a}_1 c_1 \leq\ti{a}_2 c_2 v^{-1}$. 
Therefore we have obtained characterization (i) of the proposition.\vspace*{1ex}\\
Now let $v\in\We_{N\setminus C}$ such that
\begin{eqnarray*}
 \ti{a}_1 c_1\;\leq\;\ti{a}_2 c_2 v^{-1} \quad\mb{ and }\quad \left\{\begin{array}{ccl} 
                           b_1\geq v b_2 & for & \eps=+\\
                           v^{-1} b_1\leq b_2    & for & \eps=-
\end{array}\right\}\;.
\end{eqnarray*}
Then by Lemma 2.2 of \cite{Pu3} or by Corollary \ref{P2} a) there exists an element $\ti{u}\in\We$ such that 
\begin{eqnarray*}
 \ti{a}_1\bigl(\ti{u}(c_1)^{-1}\bigr)^{-1}\;=\;\ti{a}_1 c_1 \ti{u}^{-1}\;\leq\;\ti{a}_2   \quad  \mb{ and }\quad 
    \bigl(\ti{u}(c_1)^{-1}\bigr) c_1\;=\;\ti{u}\;\leq\; c_2 v^{-1} .
\end{eqnarray*} 
Since $c_2v^{-1}\in\We_{N\setminus C}$ also $\ti{u}\in \We_{N\setminus C}$ and $\ti{u} (c_1)^{-1}\in \We_{N\setminus C}$. 
Therefore we have shown characterization (ii) of the definition of $\leq_{\eps +}$.\\
\End
Now from the last proposition and Proposition \ref{BL8} a), b), and from Remark (1) after Proposition \ref{BL3} 
follows immediately:
\begin{Prop} \label{BL11} Let $\delta\in\{+,-\}$. Let $x_1, x_2\in\We(N,C)$ and let $x_1=a_1 \idem b_1$, 
$x_2=a_2 \idem b_2$ be its normal forms II. Then $x_1\leq_{- \delta}x_2$ if and only if there exists an element 
$u\in\We_{N\setminus C}$ such that
\begin{eqnarray*}
\left\{\begin{array}{ccl}
   a_1 u^{-1}\leq a_2      & for & \delta=+\\
    a_1\geq a_2 u          & for & \delta=-
\end{array}\right\} \quad \mb{ and }\quad (x) 
\end{eqnarray*}
holds. Here for $(x)$ any of the following four statements can be taken:
\begin{eqnarray*}
\begin{array}{ccc}
(i) && b_1\leq u^{-1}b_2 \\
(ii) && u b_1\leq b_2 \vspace*{1ex}\\
(i') && b_1\leq u^{-1} \Box b_2 \\
(ii')  && u \lhd b_1\leq b_2 
\end{array}
\end{eqnarray*}
\end{Prop}
{\bf Remark:} Let ${\cal R}$ be the Renner monoid of a reductive algebraic group and let $\La\subseteq {\cal R}$ be a 
cross section lattice. Let 
$e\in\La$ and identify $\We e\We$ with $\We(N,C)$ where $N(e)=\We_N$ and $C(e)=\We_C$. Then the algebraic description 
of the Bruhat-Chevalley order obtained in \cite{PePuRe}, restricted to $\We e\We$, identifies with the 
characterization of $\leq_{++}$ of Proposition \ref{BL10}, where we take (i) for (x). From this follows by Theorem 
\ref{BL9} that the closure relation of the cells $B^\eps x B^\delta$, $x\in{\cal R}$, transfered to ${\cal R}$ and 
restricted to $\We e\We$ identifies with the extended Bruhat order $\leq_{\eps\delta}$ on $\We(N,C)$, 
$\eps,\delta\in\{+,-\}$.
In a similar way the extended Bruhat order $\leq_{\eps\delta}$ of \cite{M3}, restricted to a $\We\times\We$-orbit of 
the Weyl monoid $\widehat{W}$ identifies with the extended Bruhat order $\leq_{\eps\delta}$ here, 
$(\eps,\delta) =(+,+), (-,-),(-,+)$.
%
%
%
%
%
%
%
\section{The length of the maximal chains\label{MC}}
Immediately from Proposition \ref{BL7} follows:
\begin{Cor}\label{MC1} Let $\eps,\delta\in\{+,-\}$. Let $x,y\in\We(N,C)$ such that $x\leq_{\eps\delta}y$. The length of every 
chain joining $x$ and $y$ is finite and does not exceed $l_{\eps\delta}(y)-l_{\eps\delta}(x)$. 
In particular there exist maximal chains between $x$ and $y$.
\end{Cor}
In this section we show that every maximal $\leq_{\eps\delta}$-chain between two elements $x,y\in\We(N,C)$, 
$x\leq_{\eps\delta}y$, has length $l_{\eps\delta}(y)-l_{\eps\delta}(x)$. This also leads to the Z-Lemma for the 
extended Bruhat orders.\vspace*{1ex}\\
Usually the Bruhat order on the Coxeter group $\We$ is defined as the order is the order generated by the relations
\begin{eqnarray*}
    tx\;<\;x &\mb{ where } & x\in\We,\;t\in T\;\mb{ such that }\;l(tx)<l(x)\;.
\end{eqnarray*}
Equivalently it is the order generated by the relations
\begin{eqnarray*}
    xt\;<\;x &\mb{ where } & x\in\We,\;t\in T\;\mb{ such that }\;l(tx)<l(x)\;.
\end{eqnarray*}
The extended Bruhat orders have been defined in another way. Now we introduce a similar set of relations for the 
extended Bruhat orders. These are used for the investigation of the maximal chains of the extended Bruhat orders. 
Later we obtain as a Corollary that these relations also generate the extended Bruhat orders. 
\begin{Prop}\label{MC2} Let $\eps,\delta\in\{+,-\}$. Let $a\in \We^N$, $c\in\We_{N\setminus C}$, $b\in\mb{}^N\We$, 
and $t\in T$.
\begin{eqnarray}\label{er1}
  && \mb{For } ta<a \mb{ we have }\left\{\begin{array}{c}
    ta c\idem b <_{\eps +} ac\idem b  \\
    a c\idem b <_{\eps -} tac\idem b   \\
   \end{array}\right\}\;.  \hspace*{25ex}\\
  && \mb{For } tc<c \mb{ we have } \;a tc \idem b <_{\eps\delta} a c\idem b \;. \label{er2}  \\
  && \mb{For } bt<b \mb{ we have }\left\{\begin{array}{c}
         a c\idem b<_{+\delta} ac\idem bt  \\
         a c\idem bt<_{-\delta} ac\idem b  \\
   \end{array}\right\} \;.\label{er3}  
\end{eqnarray}
\end{Prop}
\Proof As an example we prove the first relation of (\ref{er1}). The other relations are treated in a 
similar way.
The normal form III of $tac\idem b$ is 
\begin{eqnarray*}
  (ta)^N\left((ta)_{N\setminus C}\, c\right)\idem b\;.
\end{eqnarray*}
By definition $ta c \idem b\leq_{\eps +}ac\idem b$ if there exist elements $u,v\in\We_{N\setminus C}$ such that
\begin{eqnarray*}
  (ta)^N u^{-1}\leq a\quad\mb{ and } \quad (ta)_{N\setminus C} \,c\leq u^{-1}c v^{-1}\quad\mb{ and }\quad 
\left\{\begin{array}{ccc}
  b\geq vb & \mb{for} & \eps= +\\ 
  v^{-1}b\leq b & \mb{for} & \eps= -
\end{array}\right\}\;.
\end{eqnarray*}
Since $ta\leq a$ also $(ta)^N (ta)_{N\setminus C}=(ta)^C\leq a^C=a$. Therefore the elements 
$u:=((ta)_{N\setminus C})^{-1}$ and $v:=1$ satisfy these inequalities.\\
Using Proposition \ref{BL2} we find $ta c\idem b = a c\idem b$ if and only if $a^{-1}ta\in\We_C$. Since 
$a\in\We^N\subseteq \We^C$ this would imply $ta=a\Box w$ for some $w\in\We_C$, which contradicts $ta<a$.\\
\End
\begin{Def}\label{MC3} We call a relation $x<_{\eps\delta}y$ of the form (\ref{er1}) or (\ref{er2}) 
or (\ref{er3}) elementary. We write $x<_{\eps\delta}^e y$ for short. We call a chain build by elementary 
relations an elementary chain.
\end{Def}
{\bf Remark:} For the extended Bruhat orders $\leq_{-+}$, $\leq_{++}$ it is also possible to give the elementary 
relations by using normal form I. Similar things hold for $\leq_{-+}$, $\leq_{--}$ and normal form II.\vspace*{1ex}\\
The following Theorem generalizes the Lemma of M. J. Dyer stated in 2), Section 1, to the extended Bruhat orders:
\begin{Theorem}\label{MC4} Let $\eps,\delta\in\{+,-\}$. Let $s\in S$.\\
{\bf a)} If $x<_{\eps\delta}^e y$ is an elementary relation of the form (\ref{er1}) then
\begin{eqnarray*}
   \left\{\begin{array}{ccc}
     sx<_{\eps\delta}^e sy & \mb{ if }& s\neq t \\
     sx=y                  &\mb{ if } & s=t
   \end{array}\right\} \qquad \mb{ and }\qquad xs<_{\eps\delta}^e ys\;.
\end{eqnarray*}
{\bf b)} If $x<_{\eps\delta}^e y$ is an elementary relation of the form (\ref{er2}) then
\begin{eqnarray*}
     \left\{\begin{array}{ccc}
          sx<_{\eps\delta}^e sy & \mb{ if }& s\neq ata^{-1} \\
          sx=y                  &\mb{ if } & s=ata^{-1}
     \end{array}\right\} 
 \qquad \mb{ and }\qquad
     \left\{\begin{array}{ccc}
          xs<_{\eps\delta}^e ys & \mb{ if }& s\neq (cb)^{-1}t(cb) \\
          xs=y                  &\mb{ if } & s=(cb)^{-1}t(cb)
   \end{array}\right\}  \;.
\end{eqnarray*}
{\bf c)} If $x<_{\eps\delta}^e y$ is an elementary relation of the form (\ref{er3}) then
    \begin{eqnarray*}
        sx<_{\eps\delta}^e sy
       \qquad \mb{ and }\qquad
    \left\{\begin{array}{ccc}
     xs<_{\eps\delta}^e ys & \mb{ if }& s\neq t \\
     xs=y                  &\mb{ if } & s=t
   \end{array}\right\}  \;.
\end{eqnarray*}
\end{Theorem}
\Proof We only prove the theorem for the elementary relations of the extended Bruhat order $\leq_{++}$, and for left 
multiplication by $s$. The other cases are proved in a similar way. We use several times the alternatives of V. V. 
Deohdar stated in 4), Section \ref{NP}.
Let $ac\idem b$ the normal form III of an element of $\We(N,C)$. Let $t\in T$.\vspace*{1ex}\\
{\bf To a)} Let $tac\idem b <_{++}^e ac\idem b$ be an elementary relation of the form (\ref{er1}), i.e., $ta<a$.  
If $s=t$ then trivially $s(tac\idem b)=ac\idem b$. Now let $s\neq t$. Then by 2), Section \ref{NP}, the inequality 
$ta<a$ implies
\begin{eqnarray}\label{sta}
  (sts)sa\;=\; sta\;<\;sa\;.
\end{eqnarray} 
If $sa\in\We^N$ then by (\ref{sta}) and by the definition of the elementary relations we have 
\begin{eqnarray*}
    stac\idem b \;=\; (sts)(sa)c\idem b \;<_{++}^e\;  (sa) c\idem b\;.
\end{eqnarray*}
If $sa=a\ti{s}$ with $\ti{s}\in N$ then
\begin{eqnarray}\label{stelm}
   stac\idem b =\left\{\begin{array}{lll}
            (sts)a (\ti{s}c)\idem b   & \mb{if} &  \ti{s}\in N\setminus C\\ 
            (sts)a c\idem b           & \mb{if} &  \ti{s}\in C 
   \end{array}\right\}  
\quad\mb{and}\quad  
sac\idem b=\left\{\begin{array}{lll}
   a (\ti{s}c)\idem b    & \mb{if} & \ti{s}\in N\setminus C\\
   a c\idem b            & \mb{if} & \ti{s}\in C 
    \end{array}\right\}. 
\end{eqnarray}
Since $a\in\We^N$ we have $sa=a\ti{s}=a\Box\ti{s}$. By the strong exchange condition \cite{Hu}, Section 5.8, and by the 
subword property of the Bruhat order from (\ref{sta}) follows
\begin{eqnarray*}
 sts(sa)\;=\;a \quad\mb{ or }\quad sts(sa)\;=\;a'\ti{s} \;\mb{ with }\; a'<a\;.
\end{eqnarray*}
Since $s\neq t$ it is not possible that the first equation holds. From the second we get 
\begin{eqnarray*}
   (sts)a\;=\;(sts)(sa\ti{s})\;=\; (sts (sa))\ti{s}\;=\;a'\;< \;a\;.
\end{eqnarray*} 
From this inequality and (\ref{stelm}) follows by the definition of the elementary relations 
$stac\idem b <_{++}^e sa c\idem b$.\vspace*{1ex}\\ 
{\bf To b)}  Let $atc\idem b <_{++}^e ac\idem b$ be an elementary relation of the form (\ref{er2}), i.e., $tc<c$.
Note that the subword property of the Bruhat order implies $t\in\We_{N\setminus C}$.
If $s=ata^{-1}$ then $sa tc\idem b=ac\idem b$. Now let $s\neq ata^{-1}$.\\
If $sa\in\We^N$ then by $tc<c$ and the definition of the elementary relations we get 
$ (s a)(tc)\idem b <_{++}^e (sa) c\idem b$. 
If $sa=a\ti{s}$ with $\ti{s}\in C$ then
\begin{eqnarray*} 
  satc\idem b\;=\;atc\idem b \;<_{++}^e ac\idem b\;= \;sac\idem b\;.
\end{eqnarray*}
If $sa=a\ti{s}$ with $\ti{s}\in N\setminus C$ then 
\begin{eqnarray}\label{stelm2} 
  satc\idem b\;=\;a\ti{s}tc\idem b\;=\;a\left(\ti{s}t\ti{s}(\ti{s}c)\right)\idem b \quad  \mb{ and }\quad 
  sac\idem b\;=\;a(\ti{s}c)\idem b\;. 
\end{eqnarray}
Now $s\neq ata^{-1}$ is equivalent to $\ti{s}\neq t$. By 2), Section \ref{NP}, from $tc < c$ follows 
$\ti{s}t\ti{s}(\ti{s}c)=\ti{s}tc<\ti{s}c$.  From this inequality and (\ref{stelm2}) follows by the definition 
of the elementary relations $satc\idem b <_{++}^e s a c\idem b$.\vspace*{1ex}\\ 
{\bf To c)} Let $ac\idem b <_{++}^e ac\idem bt$ be an elementary relation of the form (\ref{er3}), i.e., $bt<b$.\\
If $sa\in\We^N$ then by $bt<b$ and the definition of the elementary relations we get $sac\idem b <_{++}^e sa c\idem bt$.
If $sa=a\ti{s}$ with $\ti{s}\in C$ then 
\begin{eqnarray*}
   sac\idem b\;=\;ac\idem b\;<_{++}^e \; ac\idem bt\;=\;sac\idem bt
\end{eqnarray*}  
If $sa=a\ti{s}$ with $\ti{s}\in N\setminus C$ then
\begin{eqnarray*}
   sac\idem b\;=\;a (\ti{s}c)\idem b  \quad\mb{ and }\quad  sac\idem bt\;=\;a (\ti{s}c)\idem bt\;. 
\end{eqnarray*}
By $bt<b$ and the definition of the elementary relations we get $sac\idem b <_{++}^e sa c\idem bt$.\\
\End
Now it would be possible to prove a Z-Lemma for the orders generated by the elementary relations in the same way as 
in \cite{Dy2}, Proposition 1.9, for the twisted Bruhat orders. In our situation this is not useful because up 
to now we do not know if the elementary relations generate the extended Bruhat orders. Instead we extract in 
the next theorem certain statements about elementary chains. These are used for the inductive proof of Theorem 
\ref{MC6}, which shows the existence of a elementary chain of length $l_{\eps\delta}(y)-l_{\eps\delta}(x)$ 
between two elements $x,y\in\We(N,C)$, $x\leq_{\eps\delta}y$. 
\begin{Theorem} \label{MC5} Let $\eps,\delta\in\{+,-\}$. Let $x,y\in\We(N,C)$.\\ 
{\bf 1)} Let $s\in S$ such that $sx<_{\eps\delta} x$ and $sy<_{\eps\delta} y$. Then it is equivalent:\\
i) There exists an elementary chain of length $n$ between $x$ and $y$.\\ 
ii) There exists an elementary chain of length $n$ between $sx$ and $sy$.\vspace*{1ex}\\
{\bf 2)}  Let $s\in S$ such that $xs <_{\eps\delta} x$ and $ys <_{\eps\delta} y$. Then it is equivalent:\\
i) There exists an elementary chain of length $n$ between $x$ and $y$.\\
ii) There exists an elementary chain of length $n$ between $xs$ and $ys$.
\end{Theorem}
\Proof We prove part 1) of this theorem for the extended Bruhat order $\leq_{++}$. The other extended Bruhat 
orders are treaded similarly. Also part 2) of the theorem is proved in a similar way.\vspace*{1ex}\\
Let $ac\idem b$ an element of $\We(N,C)$ in normal form III, and let $s\in S$ such that 
\begin{eqnarray}\label{selm}
   sa c\idem b\;<_{++}\;ac\idem b\;.
\end{eqnarray} 
We show that this is already an elementary relation by using the alternatives of V. V. Deohdar which are stated 
in 4), Section \ref{NP}.
If $sa<a$ then $sa\in\We^N$. Here (\ref{selm}) is an elementary relation of the form (\ref{er1}). 
The case $sa>a$ and $sa\in\We^N$ is not possible because this would imply $sa c\idem b>_{++}ac\idem b$.
Also $sa>a$ and $sa=a\ti{s}$ with $\ti{s}\in C$ is not possible because this would imply $sa c\idem b =ac\idem b$.
If $sa>a$ and $sa=a\ti{s}$ with $\ti{s}\in N\setminus C$ then $a(\ti{s}c)\idem b$ is the normal form III of $sa c\idem b$. 
From (\ref{selm}) follows by the definition of $\leq_{++}$ that there exist elements $u,v\in\We_{N\setminus C}$ such that
\begin{eqnarray*}
  a u^{-1}\leq a \quad \mb{ and }\quad \ti{s}c \leq u^{-1}c v^{-1} \quad \mb{ and }\quad b\geq vb\;.
\end{eqnarray*}
Since $a u^{-1}= a \Box u^{-1}$ and $vb=v\Box b$, from the first and third inequality follows $u=v=1$. Inserting in 
the second inequality we get $\ti{s}c \leq c$, from which follows $\ti{s}c <c$. Therefore (\ref{selm}) is an 
elementary relation of the form (\ref{er2}).\vspace*{1ex}\\
Now we prove that $i)$ implies $ii)$. Let
\begin{eqnarray*}
   x=:z_0\;<_{++}^e\;z_1\;<_{++}^e\;\cdots\;<_{++}^e\;z_n:=y
\end{eqnarray*}
be an elementary chain of length $n$. Then
\begin{eqnarray*}
   sx=:z_{-1}\;<_{++}^e\;x=z_0\;<_{++}^e\;z_1\;<_{++}^e\;\cdots\;<_{++}^e\;z_n=y
\end{eqnarray*}
is an elementary chain of length $n+1$ with the property $sz_{-1}=z_0$. Let $p\in\{0,1,\ldots, n\}$ be the 
maximal index such that $s z_{p-1}=z_p$ resp. $z_{p-1}=s z_p$. Then due to the last theorem we get a chain of the 
form
\begin{eqnarray*}
     sx=z_{-1}\;<_{++}^e\;x=z_0\;<_{++}^e\;z_1\;<_{++}^e\cdots\;<_{++}^e\;
                                    z_{p-1}=s z_p\;<_{++}^e\;\cdots\;<_{++}^e\;sz_n=sy \;.
\end{eqnarray*}
It has length $n$.
To prove that $ii)$ implies $i)$ let 
\begin{eqnarray*}
   sx=:z_0\;<_{++}^e\;z_1\;<_{++}^e\;\cdots\;<_{++}^e\;z_n:=sy
\end{eqnarray*}
be an elementary chain of length $n$. Then
\begin{eqnarray*}
   sx=:z_0\;<_{++}^e\;<_{++}^e\;z_1\;<_{++}^e\;\cdots\;<_{++}^e\;z_n=sy\;<_{++}\;z_{n+1}:=y
\end{eqnarray*}
is an elementary chain of length $n+1$ with the property $s z_n= z_{n+1}$. Let $q\in\{0,1,\ldots, n\}$ be the 
minimal index such that $s z_q= z_{q+1}$. Then due to the last theorem we get a chain of the form
\begin{eqnarray*}
     x=sz_0\;<_{++}^e\;\cdots\;<_{++}^e\; s z_q= z_{q+1}\;<_{++}^e\;\cdots\;<_{++}^e\;z_{n+1}=y \;.
\end{eqnarray*}
It has length $n$.\\
\End
The next theorem is the key theorem of this section.
\begin{Theorem}\label{MC6}  Let $\eps,\delta\in\{+,-\}$. Let $x,y\in\We(N,C)$ such that $x\leq_{\eps\delta} y$. Then 
there exists an elementary chain between $x$ and $y$ of length $l_{\eps\delta}(y)-l_{\eps\delta}(x)$. 
\end{Theorem}
\Proof We only show this for the extended Bruhat order $\leq_{++}$, the other extended Bruhat orders are treaded 
similarly. We use several times the alternatives of V. V. Deohdar stated in 4), Section \ref{NP}.\\
Let $x=a_1c_1\idem b_1$ and $y=a_2c_2\idem b_2$ be the normal forms III of $x$ and $y$. Then by definition 
$x\leq_{++} y$ if there exist elements $u,v\in\We_{N\setminus C}$ such that 
\begin{eqnarray*}
   a_1 u^{-1}\leq a_2 \quad\mb{ and }\quad c_1\rhd v \leq u^{-1}\Box c_2  \quad\mb{ and }\quad 
   b_1\geq v b_2\;.
\end{eqnarray*}
It is easy to check from the definition of $\leq_{++}$ that 
\begin{eqnarray}\label{compch}
  a_1 c_1\idem b_1\;\leq_{++}\;a_1 c_1 v \idem b_2\;\leq_{++}\; a_1 u^{-1}c_2 \idem b_2\;\leq_{++}\;
  a_2 c_2\idem b_2\;. 
\end{eqnarray}
In the following steps 1), 2), and 3) we find elementary chains between these elements, 
whose length add up to the required length:\vspace*{1ex}\\ 
{\bf 1)} Consider in (\ref{compch}) the third inequality 
\begin{eqnarray*}
   a_1 u^{-1}c_2 \idem b_2\;\leq_{++}\;a_2 c_2\idem b_2\; \mb{ with } \;u\in\We_{N\setminus C} 
   \;\mb{ such that } \;a_1 u^{-1}\leq a_2 \;\mb{ and }\;  u^{-1}c_2=u^{-1}\Box c_2 \;.
\end{eqnarray*}
By induction over $l(a_2)$ we find an elementary chain of length $l(a_2)-l(a_1)-l(u^{-1})$ between 
$a_1 u^{-1}c_2\idem b_2$ and $a_2 c_2 \idem b_2$.\vspace*{1ex}\\
If $l(a_2)=0$ then $a_2=1$. From $1\geq a_1 u^{-1}=a_1\Box u^{-1}$ follows $a_1=1$ and $u=1$. Here $c_2\idem b_2 $ is an 
elementary chain of length $0=l(a_2)-l(a_1)-l(u^{-1})$ between $ a_1 u^{-1}c_2 \idem b_2=c_2\idem b_2 $ and 
$a_2 c_2\idem b_2=c_2\idem b_2$.\vspace*{1ex}\\
For the step of the induction choose a simple reflection $s\in S$ such that $s a_2<a_2$. In particular 
this implies $s a_2\in\We^N$.\vspace*{1ex}\\
a) Let $s a_1 < a_1$, in particular $s a_1\in\We^N$. Since $a_1 u^{-1}=a_1\Box u^{-1}$ from Theorem \ref{P2} b) 
follows $s a_1 u^{-1}< a_1 u^{-1}$. Therefore we can apply the conclusion a) of the Z-Lemma given in 1), 
Section \ref{NP}, to $a_1 u^{-1}\leq a_2$. We get $ s a_1 u^{-1}\leq s a_2$. 
By our induction assumption there exists an elementary chain of length 
$l(s a_2)-l(s a_1)-l(u^{-1})=l(a_2)-l(a_1)-l(u^{-1})$ between $s a_1 u^{-1} c_2\idem b_2$ and $s a_2 c_2\idem b_2$. Now 
by the definition of $\leq_{++}$ we have $s a_1 u^{-1} c_2\idem b_2<_{++} a_1 u^{-1} c_2\idem b_2 $ and $s a_2 c_2\idem b_2<_{++}a_2 c_2\idem b_2$. From 
Theorem \ref{MC5} follows that there also exists an elementary chain of this length between 
$a_1 u^{-1} c_2\idem b_2$ and $a_2 c_2\idem b_2$.\vspace*{1ex}\\
b) Let $s a_1 > a_1$ and $s a_1\in\We^N$. Since $s a_1 u^{-1} = s a_1 \Box u^{-1}$ from Theorem \ref{P2} b) follows 
$s a_1 u^{-1} > a_1  u^{-1}$.
Applying the conclusion c) of the Z-Lemma given in 1), Section \ref{NP}, to  $a_1 u^{-1}\leq a_2$ we get 
$a_1 u^{-1}\leq s a_2$.
By the induction assumption there exists an elementary chain of length $l(s a_2)-l(a_1)-l(u^{-1})=l(a_2)-1 -l(a_1)-l(u^{-1})$ 
between $a_1 u^{-1} c_2\idem b_2$ and $s a_2 c_2\idem b_2$. Now $s a_2 c_2\idem b_2<_{++}^e a_2 c_2\idem b_2$. By 
concatenation  there exits an elementary chain of length $l(a_2) -l(a_1)-l(u^{-1})$ between $a_1 u^{-1} c_2\idem b_2$ and 
$a_2 c_2\idem b_2$.\vspace*{1ex}\\
c) Let $s a_1 > a_1$ and $s a_1=a_1\ti{s}$ where $\ti{s}\in N$ such that $\ti{s}u^{-1}>u^{-1}$. 
Since $a_1 \ti{s} u^{-1}=a_1 \Box (\ti{s} u^{-1})$ from Theorem \ref{P2} b) follows
\begin{eqnarray*}
  s a_1 u^{-1} = a_1 \ti{s} u^{-1} > a_1  u^{-1}\;.
\end{eqnarray*} 
Now proceeding in the same way as in b) we find that there exists an elementary chain of length $l(a_2) -l(a_1)-l(u^{-1})$ 
between $a_1 u^{-1} c_2\idem b_2$ and $a_2 c_2\idem b_2$.\vspace*{1ex}\\
d) Let $s a_1 > a_1$ and $s a_1=a_1\ti{s}$ where $\ti{s}\in N$ such that $\ti{s}u^{-1}<u^{-1}$. Since 
$u\in\We_{N\setminus C}$ this last inequality is only possible for $\ti{s}\in N\setminus C$. Since 
$ a_1 u^{-1} =a_1\Box u^{-1}$ from Theorem \ref{P2} b) follows 
\begin{eqnarray*}
  a_1 u^{-1}  > a_1 \ti{s} u^{-1} = s (a_1 u^{-1}) \;.
\end{eqnarray*}
Applying conclusion a) of the Z-Lemma given in 1), Section \ref{NP}, to $a_1 u^{-1}\leq a_2$ we find
\begin{eqnarray*}
   a_1 (\ti{s} u^{-1})=s a_1 u^{-1}\leq s a_2\;.
\end{eqnarray*}
We have $\ti{s}u^{-1}c_2=\ti{s}u^{-1}\Box c_2$ and $\ti{s}u^{-1}c_2 < u^{-1}c_2$ because of
\begin{eqnarray*}
    l(\ti{s}u^{-1})+l(c_2)\;=\; -1+l(u^{-1})+l(c_2)\;=\;
   -1+ l(u^{-1}\Box c_2)\;\leq\; l(\ti{s}u^{-1}c_2)\;\leq\;l(\ti{s}u^{-1})+l(c_2)\;.
\end{eqnarray*} 
By the induction assumption there exists an elementary chain of length $l(s a_2)-l(a_1)-l(\ti{s}u^{-1})=
l(a_2)-l(a_1)-l(u^{-1})$ between $a_1 (\ti{s}u^{-1}) c_2\idem b_2$ and $s a_2 c_2\idem b_2$. Furthermore by 
the definition of $\leq_{++}$ we have
\begin{eqnarray*}
 s a_1 u^{-1}c_2\idem b_2=a_1 (\ti{s}u^{-1} c_2)\idem b_2<_{++}a_1 (u^{-1}c_2)\idem b_2  \quad \mb{ and }\quad  
 (s a_2) c_2\idem b_2<_{++} a_2 c_2\idem b_2\;.
\end{eqnarray*} 
From Theorem \ref{MC5} follows that there also exists an elementary chain of this length between 
$a_1 u^{-1} c_2\idem b_2$ and $a_2 c_2\idem b_2$.\vspace*{1ex}\\
{\bf 2)} Consider the second inequality of (\ref{compch}). Choose a maximal chain 
\begin{eqnarray*}
  c_1 v = t_m\cdots t_1 u^{-1}c_2 < \cdots< t_1 u^{-1}c_2 < u^{-1} c_2 \quad,\quad t_i\in T\cap\We_{N\setminus C}\;,
\end{eqnarray*}
in $\We_{N\setminus C}$. Then
\begin{eqnarray*}
   a_1 (c_1 v)\idem b_2  = a_1 (t_m\cdots t_1 u^{-1}c_2)\idem b_2  <_{++}^e \cdots <_{++}^e 
   a_1 (t_1 u^{-1}c_2)\idem b_2 <_{++}^e a_1 (u^{-1} c_2)\idem b_2
\end{eqnarray*}
is an elementary chain of length $m=l(u^{-1}\Box c_2)- l(c_1\rhd v)=l(u^{-1})+l(c_2)-l(c_1)+l(v)$ between 
$a_1c_1 v\idem b_2$ and $a_1u^{-1} c_2\idem b_2$.\vspace*{1ex}\\
{\bf 3)} Consider the first inequality of (\ref{compch}). Similarly as in 1), now by induction over $l(b_1)$, it is 
possible to show that there exists an elementary chain of length $l(b_1)-l(v)-l(b_2)$ between $a_1 c_1 \idem b_1$ 
and $a_1 c_1 v\idem b_2$.\vspace*{1ex}\\
By concatenation of the chains found in 1), 2), and 3) we get an elementary chain of length
\begin{eqnarray*}
   \left(l(a_2)-l(a_1)-l(u^{-1})\right)+\left(l(u^{-1})+l(c_2)-l(c_1)+l(v)\right)+\left(l(b_1)-l(v)-l(b_2)\right)\\
  \;=\;  \left(l(a_2)+l(c_2)-l(b_2) \right) - \left(l(a_1)+l(c_1)-l(b_1)\right)\;=\; 
   l_{++}(a_2 c_2\idem b_2)-l_{++}(a_1 c_1\idem b_1)
\end{eqnarray*}
between $a_1 c_1\idem b_1$ and $a_2 c_2\idem b_2$.\\
\End
As an immediate consequence of Theorem \ref{MC6} we get:
\begin{Cor}\label{MC8}  Let $\eps,\delta\in\{+,-\}$. The order relation $\leq_{\eps\delta}$ is generated by its 
elementary relations. 
\end{Cor}
Also as a consequence of Theorem \ref{MC6} we obtain the lengths of the maximal chains of the extended Bruhat orders 
between two elements:
\begin{Cor}\label{MC7} Let $\eps,\delta\in\{+,-\}$. Let $x,y\in\We(N,C)$ with $x\leq_{\eps\delta} y$. Every 
maximal chain between $x$ and $y$ is elementary and has length $l_{\eps\delta}(y)-l_{\eps\delta}(x)$.
\end{Cor}
\Proof Since the order relation $\leq_{\eps\delta}$ is generated by its 
elementary relations, every maximal chain between $x$ and $y$ is elementary. 
Now let 
\begin{eqnarray*}
  x=z_0 <_{\eps\delta} z_1 <_{\eps\delta}\cdots <_{\eps\delta} z_m=y 
\end{eqnarray*}
be a maximal chain. By Proposition \ref{BL7} we have $l_{\eps\delta}(z_{i-1})<l_{\eps\delta}(z_i)$ for all 
$i\in\{1,2,\ldots ,m\}$. Now suppose that there exists an index $i$ such that 
$l_{\eps\delta}(z_{i-1})+2\leq l_{\eps\delta}(z_i)$. Then by Theorem \ref{MC6} there exists an element 
$z\in\We(N,C)$ such that $ z_{i-1}<_{\eps\delta} z <_{\eps\delta} z_i$, which contradicts the maximality of the 
chain. 
Therefore $l_{\eps\delta}(z_{i-1})+1=l_{\eps\delta}(z_i)$ for all $i\in\{1,2,\ldots ,m\}$, from which follows 
$m=l_{\eps\delta}(y)-l_{\eps\delta}(x)$. \\
\End
In the same way as in \cite{Dy2}, Proposition 1.9, for the twisted Bruhat orders, it is now possible to complete the 
proof of the Z-Lemma for the extended Bruhat orders:
\begin{Cor}[Z-Lemma]  Let $\eps,\delta\in\{+,-\}$. Let $x,y\in\We(N,C)$.\vspace*{1ex}\\
{\bf 1)} For $s\in S$ such that $sx<_{\eps\delta} x$ and $sy<_{\eps\delta} y$ the following statements are 
equivalent:\\
i) $x\leq_{\eps\delta} y$. \\
ii) $sx\leq_{\eps\delta} y$. \\
iii) $sx\leq_{\eps\delta} sy$. \\
{\bf 2)} For $s\in S$ such that $xs<_{\eps\delta} x$ and $ys<_{\eps\delta} y$ the following statements are 
equivalent:\\
i) $x\leq_{\eps\delta} y$. \\
ii) $xs\leq_{\eps\delta} y$. \\
iii) $xs\leq_{\eps\delta} ys$. 
\end{Cor}
%
%
%
%
%
%
%
%
%
%
%
%
%
%
%
%
%
%

%
%
%
%

\begin{thebibliography}{AAAA1}
\addcontentsline{toc}{section}{References}
%
%
\bibitem[Ch,Dy]{ChDy}{Y. Chen, M. J. Dyer, On the combinatorics of $B\times B$-orbits on group compactifications, Journal of Algebra {\bf 263} (2003), 278-293}
%
%
\bibitem[De]{De}{V. V. Deodhar, Some Characterizations of Bruhat Ordering on a Coxeter group and Determination of the 
Relative M\"obius Function, Inventiones math., {\bf 39} (1977), 187-198}
%
\bibitem[Dy 1]{Dy1}{M. J. Dyer, Iwahori Hecke Algebras and Shellings of Bruhat Intervals, 
Compositio Math. {\bf 89} (1993), 91-115}
%
%
\bibitem[Dy 2]{Dy2}{M. J. Dyer, Hecke Algebras and Shellings of Bruhat Intervals II; Twisted Bruhat Orders, 
Contemporary Mathematics, {\bf 139} (1992), 141-165}
%
\bibitem[Dy 3]{Dy3}{M. J. Dyer, Quotients of Twisted Bruhat Orders, Journal of Algebra {\bf 163} (1994), 861-879}
%
\bibitem[Hu]{Hu}{J. E. Humphreys, Reflection Groups and Coxeter Groups, Cambridge Studies in Advanced Mathematics 29, Cambridge University 
Press, 1990}
%
%
\bibitem[K,P]{KP}{V. G. Kac, D. H. Peterson, Regular functions on certain in\-finite-\-di\-men\-si\-onal groups, in Arithmetic and Geometry, 
Progress in Math. {\bf 36}, Birkh\"auser, Boston, 1983, 141-166}
%
%
%
\bibitem[Kas]{Kas}{M. Kashiwara, The Flag Manifold of Kac-Moody Lie Algebra,
American Journal of Mathematics {\bf 111} (1989), Supplement Volume: Algebraic 
Analysis, Geometry and Number Theory, Ed J. i. Jgusa, John Hopkins University, 
Baltimore, 1989, 161-190}
%
%
%
\bibitem[M 1]{M1}{C. Mokler, An analogue of a reductive algebraic monoid, whose unit group is a Kac-Moody group, 
arXiv.org e-print: math.AG/0204246, submitted to Memoirs of the AMS}
%
\bibitem[M 2]{M2}{C. Mokler, The $\F$-valued points of the algebra of strongly regular functions of a Kac-Moody group, 
Transformation Groups {\bf 7} (2002), 343-378}
%
\bibitem[M 3]{M3}{C. Mokler, Extending the Bruhat order and the length function from the Weyl group to the Weyl monoid
arXiv.org e-print: math.RT/0303281, accepted by the Journal of Algebra}
%
%
%
%
\bibitem[Pi]{Pi}{D. Pickrell, Invariant Measures for Unitary Groups 
associated to Kac-Moody Lie Algebras, Memoirs of the AMS {\bf 146} (2000)}
%
%
%
\bibitem[Pe,Pu,Re]{PePuRe}{E. A. Pennel, M. S. Putcha, L. E. Renner, Analogue of the Bruhat-Chevalley order for 
reductive monoids, Journal of Algebra, {\bf 196} (1997), 339-368}
%
%
\bibitem[Pu 1]{Pu1}{M. Putcha, Idempotent cross-sections of J-classes, Semigroup Forum {\bf 26} (1983), 
103-109}
%
\bibitem[Pu 2]{Pu2}{M. Putcha, Reductive groups and regular semigroups, Semigroup Forum {\bf 30} (1984), 253-261}
%
%
\bibitem[Pu 3]{Pu3}{M. S. Putcha, Shellability in reductive monoids, Transactions of the American Mathematical 
Society, {\bf 354} (2002), 413-426}
%
\bibitem[Pu, Re]{PuRe}{M. S. Putcha, L. E. Renner, The Systems of Idempotents and the Lattice of 
J-classes of Reductive Algebraic Monoids, Journal of Algebra {\bf 116} (1988), 385-399}
%
\bibitem[Re 1]{Re1}{L. E. Renner, Analogue of the Bruhat Decomposition for Algebraic Monoids, Journal of 
Algebra {\bf 101} (1986), 303-338}
%
\bibitem[Re 2]{Re2}{L. E. Renner, Analogue of the Bruhat Decomposition for Algebraic Monoids II, The length 
function and trichotomy, Journal of Algebra {\bf 175} (1995), 695-714}
%
\bibitem[Re 3]{Re3}{L. E. Renner, An explicit Cell Decomposition of the Wonderful Compactification of a Semisimple 
Algebraic Group, Canadian Mathematical Bulletin {\bf 46} (2003),140-148}
%
%
%
\bibitem[So]{So}{L. Solomon, The Bruhat decomposition, Tits system and Iwahori ring for the monoid of matrices over a finite field, Geom. Dedicata {\bf 30} (1990), 15-49}
%
%
\bibitem[Sp]{Sp}{T. A. Springer, Intersection cohomology of $B\times B$-orbit closures in group compactifications, Journal of Algebra {\bf 258} (2002), 71-111}
%
%
%
\end{thebibliography}
\end{document}